\documentclass[11pt,a4paper]{article}
 \usepackage{emlines2,amstext,amssymb,amscd,amsmath,color}

\def\Pin{P^{\rm in}}
\def\Pout{P^{\rm out}}
\def\ellin{t}
\def\ellout{h}

\def\sprod{\;\hbox{\unitlength=1mm\begin{picture}(4,4)%
\put(0,0){$\triangleright$}\put(1.5,0){$\triangleleft$}%
\end{picture}}\;}

\begin{document}
\title{On the structure of regular
$B_2$-type crystals} \author{ V.I.Danilov\thanks{Central Institute
of Economics and Mathematics, Russian Academy of Sciences,
Nakhimovskij Prospekt 47, Moscow 117418, Russia. Email:
danilov@cemi.rssi.ru}, A.V.Karzanov\thanks{Institute for System
Analysis of the RAS, 9, Prospect 60 Let Oktyabrya, 117312 Moscow,
Russia. Email: sasha@cs.isa.ru. }, and G.A.Koshevoy\thanks{Central
Institute of Economics and Mathematics, Russian Academy of
Sciences, and Poncelet laboratory (UMI 2615 du CNRS). Email:
koshevoy@cemi.rssi.ru}}

\date{}
\maketitle

\parskip=3pt

\section{Introduction}

For simply-laced Kac-Moody algebras $\frak g$, Stembridge \cite{stem} proposed a `local' axiomatization of
crystal graphs of representations of $U_q(\frak g)$. In fact, because of an important result in \cite{KKM-92} an
essential part in studying such crystals should have been carried out for the simplest nontrivial case
$\frak g=sl(3)$. Our paper \cite{a2} gives a combinatorial construction that describes the formation of any
crystal of representations of $U_q(sl(3))$, a regular crystal graph of type $A_2$.

In this paper we attempt to carry out a similar programme for the
algebra  $sp(4)$.
We follow the ideology of \cite{a2} and propose axioms of
monotonicity and commutativity for edge-2-colored graphs which
characterize the crystals of integrable representations of
$U_q(sp(4))$, regular crystal graphs of $B_2$-type. An
edge-colored directed graph which obeys our Axioms (K0)--(K5) is
called an R-{\em graph} (for brevity), and our main result is that
the regular crystals of $B_2$-type are R-graphs and vice versa. In
particular, axiom (K5) refines Sternberg's $B_2$-type relations by
involving a certain labeling on the crystal edges. At the end of
\cite{stem}, Stembridge conjectured a list of relations between
crystal operations for the $B_2$ case and Sternberg proved this
conjecture in \cite{stern}.

Further, we give a direct combinatorial construction for the
crystals in question, using a general operation on graphs from
\cite{a2}. On this way we introduce a new, so-called {\em crossing
model}, which does not exploit Young tableaux. This combinatorial
model consists of a two-component graph of a rather simple form
and of a certain set of integer-valued functions on its vertices,
and we show that these functions one-to-one correspond to the
vertices of a regular $B_2$-crystal. In particular, due to this
model, the vertices of a regular $B_2$-crystal can be located at
the integer points of the union of five 4-dimensional polyhedra in
$\mathbb R^7$, though this union is not a convex set. (Another
crossing model is developed for regular $A_n$-crystals
in~\cite{DKK-06a}.)

The structure of this paper is the following. Section 2 is devoted
to axioms. In Section 3 we formulate the main result and give  a
constructive characterization of the regular $B_2$-crystal.
In Section 4 we introduce the crossing model for such crystals.
In Section 5 we prove that certain
intervals of the admissible configurations on this model are just
regular $B_2$-crystals, relying on the known
characterization of the latter via Littelmann's cones. In Section 6
we prove
that the graphs generated by the intervals on the crossing model
are essentially the same as those constructed in Section 3.

{\bf Acknowledgements.}  The authors gratefully acknowledge the
support of the grant 047.011.2004.017 from NWO and RFBR, and the
grant 05-01-02805 from CNRSLa from CNRS and RFBR. G.~Koshevoy
gratefully acknowledges the support of Support Foundation of the
Russian Federation and IHES for hospitality and perfect working
conditions. Part of this research was done while A.~Karzanov was
visiting PNA1 at CWI, Amsterdam and was supported by a grant from this Center.

\section{Axiomatization of $B_2$-type regular crystals}

Here we consider edge-two-colored directed graphs managed by the
$B_2$-type Cartan matrix $A=\begin{pmatrix} 2& -2\cr -1
&2\cr\end{pmatrix}$. Of our interest will be those graphs that
match the list of axioms below.

\medskip
We consider a digraph $G=(V,A_1\coprod A_2)$ with the vertex set
$V$ and the edges set partitioned into two subsets $A_1$ and
$A_2$. For convenience we refer to the edges in $A_1$ as being
colored in color 1, say {\em red}, and the edges in $A_2$ as being
colored in color 2, say {\em green}.
\medskip

The first axiom concerns the structure of the monochromatic graphs
$(V,A_1 )$ and $(V,A_2 )$, and states that they are constituted as
the disjoint union of finite monochromatic strings. Specifically,
\medskip

\begin{itemize}
\item[(K0)] For $i=1,2$, each (weakly) connected component of
$(V,A_i)$ is a finite simple (directed) path, i.e., a sequence of
the form $(v_0,e_1,v_1,\ldots,e_k,v_k)$, where
$v_0,v_1,\ldots,v_k$ are distinct vertices and each $e_i$ is an
edge going from $v_{i-1}$ to $v_i$.
  \end{itemize}

In particular, each vertex $v$ has at most one outgoing 1-colored
edge and at most one incoming 1-colored edge, and similarly for
2-colored edges. For brevity, we refer to a {\em maximal}
monochromatic path in $K$, with color $i$ on the edges, as an
$i$-{\em string}. The $i$-string passing through a given vertex
$v$ (possibly consisting of the only vertex $v$) is denoted by
$P_i(v)$, its part from the first vertex to $v$ by $\Pin_i(v)$,
and its part from $v$ to the last vertex by $\Pout_i(v)$. The
lengths of $\Pin_i(v)$ and of $\Pout_i(v)$ (i.e., the numbers of
edges in these paths) are denoted by $\ellin_i(v)$ and
$\ellout_i(v)$, respectively. The next axiom tells how these
lengths can change when one traverses an edge of the other color.

Consider a red edge $(x,y)$ connecting a vertex $x$ to a vertex
$y$ (using the operator notation for Kashiwara crystals, one can
write $y=F_1x$), and consider the strings $P_2(x)$ and $P_2(y)$.
For a
2-colored edge $(p,q)$ going from $p$ to $q$, $q=F_2p$, we
consider the pair of strings $P_1(p)$ and $P_1(q)$. 
\begin{itemize}
\item[(K1)]
a) If $y=F_1(x)$, then $t_2(x)-t_2(y)+h_2(y)-h_2(x)=1$.

b) If $q=F_2(p)$, then $t_1(p)-t_1(q)+h_1(q)-h_1(p)=2$.

\end{itemize}

Let us define the weight function $wt:V\to\mathbb Z^2$ by the rule
$$
w\colon x\to (h_1(x)-t_1(x), h_2(x)- t_2(x)).
$$
Then, the functions $t_i$, $h_i$, $i=1,2$, and $wt$ define a
semi-normal crystal of $B_2$-type due to Kashiwara \cite{kas-95}. \medskip

The next axioms is due to Stembridge \cite{stem}.

\begin{itemize}
\item[(K2)] For each 1-colored edge $e$ going from $x$ to $y$,
$y=F_1(x)$, there hold $t_2(x)-t_2(y)\ge 0$ and
${h_2(y)-h_2(x)}\ge 0$.

For each 2-colored edge $e'$ going from $p$ to $q$, $q=F_2p$, and
there hold $t_1(p)-t_1(q)\ge 0$ and ${h_1(q)-h_1(p)}\ge 0$.
\end{itemize}

By this reason, we endow the edges of $G$ with {\em labels}, 
and for this reason, we say that $G$ is a {\em
decorated graph}. More precisely, the red edges have 
labels $0$ or $1$ and the label being assigned to a red
edge $e=(x,y=F_1(x))$ is equal to $h_2(y)-h_2(x)$. 
The green edges have labels $0$ or $\frac 12$ or $1$ and the
label assigned to a green edge $e'=(p,q)$ is equal 
to $\frac{h_1(q)-h_1(p)}2\in\{0,\frac 12, 1\}$.

The following axioms were found experimentally and their
justification is done by Theorem 1 in which we assert that regular
crystal graphs of $B_2$-type characterized by the whole list of
axioms.

The third axiom is

\begin{itemize}
\item[(K3)] For any $i$ and any $i$-string $P_i$, the labels on consecutive
edges $e_1$ and $e_2$ (where $e_2$ follows $e_1$) do not decrease.
\end{itemize}

The next axiom presents  important commutation between special
pairs of green and red edges.

\begin{itemize}
\item[(K4)] Suppose a red edge enters a vertex $v$ and has label $a$ and let
a green edge leave $v$ and have label $b$. Then there holds $b\neq
a$, moreover, if $ a< b$, then $a=0$, $b=1$ and the commutative
diagram takes place:
\bigskip

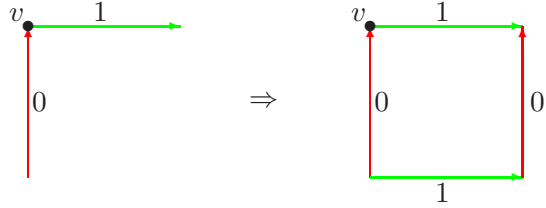
\unitlength=0.5mm \special{em:linewidth 0.4pt}
\linethickness{0.4pt}
\begin{picture}(84.00,54.00)(0,30)
\put(110.00,40.00){\color{red}\vector(0,1){40.00}}
\put(110.00,80.00){\color{green}\vector(1,0){40.00}}
\put(110.00,40.00){\color{green}\vector(1,0){40.00}}
\put(150.00,40.00){\color{red}\vector(0,1){40.00}}
\put(107.00,83.00){\makebox(0,0)[cc]{$v$}}
\put(113.00,60.00){\makebox(0,0)[cc]{$0$}}
\put(129.00,84.00){\makebox(0,0)[cc]{$1$}}
\put(154.00,60.00){\makebox(0,0)[cc]{$0$}}
\put(129.00,36.00){\makebox(0,0)[cc]{$1$}}
\put(110.00,80.00){\circle*{2.8}}
\put(20.00,40.00){\color{red}\vector(0,1){40.00}}
\put(20.00,80.00){\color{green}\vector(1,0){40.00}}
\put(20.00,80.00){\circle*{2.8}}
\put(17.00,83.00){\makebox(0,0)[cc]{$v$}}
\put(23.00,60.00){\makebox(0,0)[cc]{$0$}}
\put(39.00,84.00){\makebox(0,0)[cc]{$1$}}
\put(82.00,60.00){\makebox(0,0)[cc]{$\Rightarrow$}}
\end{picture}

Similarly, suppose a green edge enters a vertex $q$ and has label
$b$ and a red edge leaves $q$ and has label $a$. Then there holds
$b\neq a$, and if $ a> b$, then $a=1$, $b=0$, and the commutative
diagram takes place:\bigskip

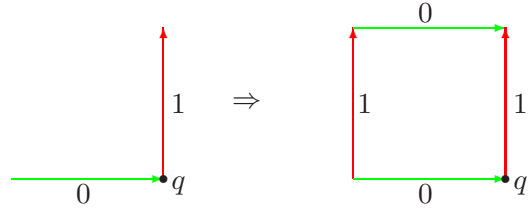
\unitlength=0.5mm \special{em:linewidth 0.4pt}
\linethickness{0.4pt}
\begin{picture}(84.00,54.00)(0,30)
\put(110.00,40.00){\color{red}\vector(0,1){40.00}}
\put(110.00,80.00){\color{green}\vector(1,0){40.00}}
\put(110.00,40.00){\color{green}\vector(1,0){40.00}}
\put(150.00,40.00){\color{red}\vector(0,1){40.00}}
\put(113.00,60.00){\makebox(0,0)[cc]{$1$}}
\put(129.00,84.00){\makebox(0,0)[cc]{$0$}}
\put(154.00,60.00){\makebox(0,0)[cc]{$1$}}
\put(129.00,36.00){\makebox(0,0)[cc]{$0$}}
\put(150.00,40.00){\circle*{2.00}}
\put(154.00,38.00){\makebox(0,0)[cc]{$q$}}
\put(20.00,40.00){\color{green}\vector(1,0){40.00}}
\put(60.00,40.00){\color{red}\vector(0,1){40.00}}
\put(60.00,40.00){\circle*{2.00}}
\put(64.00,38.00){\makebox(0,0)[cc]{$q$}}
\put(39.00,36.00){\makebox(0,0)[cc]{$0$}}
\put(64.00,60.00){\makebox(0,0)[cc]{$1$}}
\put(82.00,60.00){\makebox(0,0)[cc]{$\Rightarrow$}}
\end{picture}

\end{itemize}

From these axioms we have the following property of graphs which satisfy
axioms $K1-K4$:

A green edge labeled 1 is at the top of the corresponding
commutative square in Axiom 4; a green edge labeled 0 is at the
bottom the corresponding commutative square in Axiom 4. A green
edge labeled $\frac 12$ is always a part of the crystal graph of a
fundamental representation of $U_q(sp(4))$.

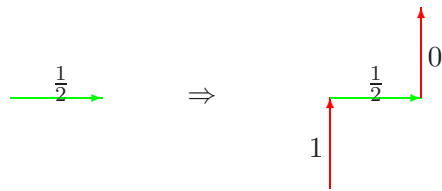
\unitlength=0.6mm \special{em:linewidth 0.4pt}
\linethickness{0.4pt}
\begin{picture}(113.00,60.00)(0,50)
\put(20.00,80.00){\color{green}\vector(1,0){20.00}}
\put(31.00,83.00){\makebox(0,0)[cc]{$\frac 12$}}
\put(90.00,60.00){\color{red}\vector(0,1){20.00}}
\put(90.00,80.00){\color{green}\vector(1,0){20.00}}
\put(110.00,80.00){\color{red}\vector(0,1){20.00}}
\put(87.00,69.00){\makebox(0,0)[cc]{$1$}}
\put(100.00,83.00){\makebox(0,0)[cc]{$\frac 12$}}
\put(113.00,89.00){\makebox(0,0)[cc]{$0$}}
\put(62.00,80.00){\makebox(0,0)[cc]{$\Rightarrow$}}
\end{picture}

In fact, because of axiom (K2), there is an ingoing red edge to
the starting point of a green edge labeled $\frac 12$, and there
is an outgoing red edges from the ending point of the green edge.
Due to axiom K4, the ingoing red edge can not be labeled 0 and the
outgoing can not be labeled 1. Thus they are labeled 1 and 0
respectively.\medskip

\noindent{\bf Corollary}. {\em In an S-graph any green string has at most one edge labeled $\frac 12$}.

{\em Proof}.  Suppose two green edges, both labeled $\frac 12$,
belong to a green string. Then, because of Axiom K3, these edges
have  a common vertex. Then due to above property of green edges
labeled $\frac 12$, the ingoing red edge to this vertex is labeled
1 and the outgoing edge is labeled 0. That contradicts to the
monotonicity Axiom K3.\medskip

Because to Axiom K3 each red string has a unique vertex (possibly
the beginning or ending one) at which labels switch from $0$ to
$1$, we call this vertex {\em critical}. Because of Axiom 3 and
Corollary each green string has either a critical vertex or a {\em
critical edge} labeled $\frac 12$.
\medskip


Because of the above property of S-graphs, of interest are the
green edges labeled 1 which are not bottom edges of commutative
squares (of course, each such an edge is at the top a commutative
square) and the green edges labeled 0 which are not top edges of
the corresponding commutative squares. Possible configurations
"above" and "below" such edges, respectively, are managed by the
following axiom.

\begin{itemize}
\item[(K5)] The implications illustrated on  Pictures 1,
2 and dual to these implications have to hold.
\end{itemize}

On Pictures 1 and 2 we demonstrate four possible configurations
which involve a green edge labeled 1, which does not located at
the bottom of a commutative square. Possible configurations
involving a green edge labeled 0, which does not located at the
top of a commutative square, are obtained dually to that indicated
at Pictures 1 and 2, that is one has to reverse directions of
edges and replace a label $a$ to $1-a$.
\medskip

\unitlength=1.00mm \special{em:linewidth 0.4pt}
\linethickness{0.4pt}
\begin{picture}(124.00,153.00)
\put(20.00,120.00){\color{green}\vector(1,0){15.00}}
\put(20.00,120.00){\color{red}\vector(0,1){15.00}}
\put(20.00,135.00){\color{green}\vector(1,0){13.00}}
\put(26.00,117.00){\makebox(0,0)[cc]{$1$}}
\put(17.00,127.00){\makebox(0,0)[cc]{$1$}}
\put(26.00,138.00){\makebox(0,0)[cc]{$0$}}
\put(57.00,126.00){\makebox(0,0)[cc]{$\Rightarrow$}}
\put(20.00,80.00){\color{green}\vector(1,0){15.00}}
\put(35.00,80.00){\color{red}\vector(0,1){15.00}}
\put(22.00,95.00){\color{green}\vector(1,0){13.00}}
\put(27.00,76.00){\makebox(0,0)[cc]{$1$}}
\put(38.00,86.00){\makebox(0,0)[cc]{$0$}}
\put(28.00,98.00){\makebox(0,0)[cc]{$\frac 12$}}
\put(57.00,87.00){\makebox(0,0)[cc]{$\Rightarrow$}}
\put(20.00,30.00){\color{green}\vector(1,0){15.00}}
\put(35.00,30.00){\color{red}\vector(0,1){15.00}}
\put(35.00,45.00){\color{red}\vector(0,1){15.00}}
\put(20.00,30.00){\color{red}\vector(0,1){17.00}}
\put(20.00,47.00){\color{green}\vector(1,0){17.00}}
\put(37.00,47.00){\circle{1.50}}
\put(27.00,26.00){\makebox(0,0)[cc]{$1$}}
\put(17.00,37.00){\makebox(0,0)[cc]{$1$}}
\put(37.00,37.00){\makebox(0,0)[cc]{$0$}}
\put(28.00,50.00){\makebox(0,0)[cc]{$\frac 12$}}
\put(57.00,42.00){\makebox(0,0)[cc]{$\Rightarrow$}}
\put(75.00,104.00){\makebox(0,0)[cc]{a)}}
\put(75.00,70.00){\makebox(0,0)[cc]{b)}}
\put(75.00,5.00){\makebox(0,0)[cc]{c)}}
\put(90.00,120.00){\color{green}\vector(1,0){15.00}}
\put(90.00,120.00){\color{red}\vector(0,1){15.00}}
\put(90.00,135.00){\color{green}\vector(1,0){13.00}}
\put(96.00,117.00){\makebox(0,0)[cc]{$1$}}
\put(87.00,127.00){\makebox(0,0)[cc]{$1$}}
\put(96.00,138.00){\makebox(0,0)[cc]{$0$}}
\put(104.00,135.00){\color{green}\vector(1,0){16.00}}
\put(120.00,135.00){\color{red}\vector(0,1){15.00}}
\put(105.00,120.00){\color{red}\vector(0,1){13.00}}
\put(105.00,133.00){\color{red}\vector(0,1){17.00}}
\put(105.00,150.00){\color{green}\vector(1,0){15.00}}
\put(105.00,133.00){\circle{1.50}}
\put(113.00,153.00){\makebox(0,0)[cc]{$0$}}
\put(122.00,143.00){\makebox(0,0)[cc]{$0$}}
\put(103.00,144.00){\makebox(0,0)[cc]{$1$}}
\put(113.00,133.00){\makebox(0,0)[cc]{$1$}}
\put(107.00,127.00){\makebox(0,0)[cc]{$0$}}
\put(106.00,95.00){\color{green}\vector(1,0){15.00}}
\put(121.00,95.00){\color{red}\vector(0,1){15.00}}
\put(107.00,110.00){\color{green}\vector(1,0){14.00}}
\put(113.00,91.00){\makebox(0,0)[cc]{$1$}}
\put(124.00,101.00){\makebox(0,0)[cc]{$0$}}
\put(114.00,113.00){\makebox(0,0)[cc]{$\frac 12$}}
\put(107.00,94.00){\circle{1.50}}
\put(107.00,80.00){\color{red}\vector(0,1){14.00}}
\put(107.00,94.00){\color{red}\vector(0,1){16.00}}
\put(105.00,95.00){\circle*{2.00}}
\put(90.00,95.00){\color{green}\vector(1,0){15.00}}
\put(90.00,80.00){\color{green}\vector(1,0){17.00}}
\put(90.00,80.00){\color{red}\vector(0,1){15.00}}
\put(96.00,97.00){\makebox(0,0)[cc]{$0$}}
\put(86.00,88.00){\makebox(0,0)[cc]{$1$}}
\put(97.00,76.00){\makebox(0,0)[cc]{$\frac 12$}}
\put(105.00,87.00){\makebox(0,0)[cc]{$0$}}
\put(105.00,102.00){\makebox(0,0)[cc]{$1$}}
\put(85.00,30.00){\color{green}\vector(1,0){15.00}}
\put(100.00,30.00){\color{red}\vector(0,1){15.00}}
\put(100.00,45.00){\color{red}\vector(0,1){15.00}}
\put(85.00,30.00){\color{red}\vector(0,1){17.00}}
\put(85.00,47.00){\color{green}\vector(1,0){17.00}}
\put(103.00,47.00){\circle{1.50}}
\put(92.00,26.00){\makebox(0,0)[cc]{$1$}}
\put(82.00,37.00){\makebox(0,0)[cc]{$1$}}
\put(102.00,37.00){\makebox(0,0)[cc]{$0$}}
\put(93.00,50.00){\makebox(0,0)[cc]{$\frac 12$}}
\put(88.00,45.00){\color{green}\vector(1,0){12.00}}
\put(88.00,45.00){\color{red}\vector(0,1){13.00}}
\put(103.00,47.00){\color{red}\vector(0,1){11.00}}
\put(88.00,58.00){\color{green}\vector(1,0){15.00}}
\put(20.00,30.00){\circle*{2.00}}
\put(85.00,30.00){\circle*{2.00}}
\put(93.00,42.00){\makebox(0,0)[cc]{$\frac 12$}}
\put(94.00,60.00){\makebox(0,0)[cc]{$0$}}
\put(105.00,52.00){\makebox(0,0)[cc]{$0$}}
\put(34.00,135.00){\circle{1.50}}
\put(104.00,135.00){\circle{1.50}}
\put(20.00,120.00){\circle*{2.00}}
\put(36.00,120.00){\circle*{2.00}}
\put(90.00,120.00){\circle*{2.00}}
\put(106.00,120.00){\circle*{2.00}}
\put(19.00,80.00){\circle*{2.00}}
\put(36.00,80.00){\circle*{2.00}}
\put(122.00,95.00){\circle*{2.00}}
\put(36.00,30.00){\circle*{2.00}}
\put(101.00,30.00){\circle*{2.00}}
\put(88.00,45.00){\circle{1.50}}
\put(86.00,52.00){\makebox(0,0)[cc]{$1$}}
\put(88.00,32.00){\color{red}\vector(0,1){13.00}}
\put(88.00,33.00){\circle{2.00}}
\put(90.00,38.00){\makebox(0,0)[cc]{$1$}}
\put(33.00,54.00){\makebox(0,0)[cc]{$0$}}
\put(98.00,54.00){\makebox(0,0)[cc]{$0$}}
\end{picture}
\begin{center}
Picture 1.
\end{center}

\unitlength=0.8mm \special{em:linewidth 0.4pt}
\linethickness{0.4pt}
\begin{picture}(128.00,75.00) \put(78.00,11.00){\circle*{2.00}}
\put(78.00,11.00){\color{green}\vector(1,0){15.00}}
\put(93.00,11.00){\color{red}\vector(0,1){15.00}}
\put(93.00,26.00){\color{red}\vector(0,1){15.00}}
\put(93.00,41.00){\color{green}\vector(1,0){15.00}}
\put(108.00,41.00){\color{red}\vector(0,1){17.00}}
\put(108.00,58.00){\color{green}\vector(1,0){17.00}}
\put(125.00,58.00){\color{red}\vector(0,1){16.00}}
\put(78.00,11.00){\color{red}\vector(0,1){17.00}}
\put(78.00,28.00){\color{green}\vector(1,0){13.00}}
\put(91.00,28.00){\color{red}\vector(0,1){16.00}}
\put(91.00,44.00){\color{green}\vector(1,0){15.00}}
\put(106.00,44.00){\color{red}\vector(0,1){16.00}}
\put(106.00,60.00){\color{red}\vector(0,1){14.00}}
\put(106.00,74.00){\color{green}\vector(1,0){19.00}}
\put(108.00,41.00){\circle*{0.00}}
\put(91.00,44.00){\circle*{2.00}}
\put(125.00,74.00){\circle*{2.00}}
\put(75.00,19.00){\makebox(0,0)[cc]{$1$}}
\put(84.00,31.00){\makebox(0,0)[cc]{$\frac 12$}}
\put(89.00,36.00){\makebox(0,0)[cc]{$0$}}
\put(98.00,46.00){\makebox(0,0)[cc]{$1$}}
\put(104.00,52.00){\makebox(0,0)[cc]{$0$}}
\put(104.00,67.00){\makebox(0,0)[cc]{$1$}}
\put(115.00,77.00){\makebox(0,0)[cc]{$0$}}
\put(128.00,67.00){\makebox(0,0)[cc]{$0$}}
\put(117.00,55.00){\makebox(0,0)[cc]{$\frac 12$}}
\put(110.00,49.00){\makebox(0,0)[cc]{$1$}}
\put(102.00,37.00){\makebox(0,0)[cc]{$0$}}
\put(95.00,32.00){\makebox(0,0)[cc]{$1$}}
\put(95.00,18.00){\makebox(0,0)[cc]{$0$}}
\put(86.00,7.00){\makebox(0,0)[cc]{$1$}}
\put(15.00,11.00){\circle*{2.00}}
\put(15.00,11.00){\color{green}\vector(1,0){15.00}}
\put(30.00,11.00){\color{red}\vector(0,1){15.00}}
\put(30.00,26.00){\color{red}\vector(0,1){15.00}}
\put(15.00,11.00){\color{red}\vector(0,1){17.00}}
\put(15.00,28.00){\color{green}\vector(1,0){13.00}}
\put(12.00,19.00){\makebox(0,0)[cc]{$1$}}
\put(21.00,31.00){\makebox(0,0)[cc]{$\frac 12$}}
\put(32.00,32.00){\makebox(0,0)[cc]{$1$}}
\put(32.00,18.00){\makebox(0,0)[cc]{$0$}}
\put(23.00,7.00){\makebox(0,0)[cc]{$1$}}
\put(74.00,8.00){\makebox(0,0)[cc]{$A$}}
\put(90.00,47.00){\makebox(0,0)[cc]{$B$}}
\put(127.00,77.00){\makebox(0,0)[cc]{$C$}}
\put(113.00,42.00){\makebox(0,0)[cc]{$D$}} 
\put(108.00,41.00){\circle*{2.00}}
\put(63.00,44.00){\makebox(0,0)[cc]{$\Rightarrow$}}
\end{picture}

\begin{center} Picture 2.
\end{center}

{\bf Remark}. The relation on Picture 2 is the Verma $B_2$-type
relation of degree $7$ (if we ignore labels), though depicted in
the form somewhat different from that of the custom Verma
$B_2$-type relation. This form of the Verma relation can be
interpreted as commuting of crystals of the fundamental
representations of $U_q(sp(4))$, the paths between $A$ and $B$,
and $B$ and $C$, and the paths between $A$ and $D$, and $D$ and $C$, respectively.

The custom Verma relation is obtained if one draws two commutative
diagrams which follow from Axiom K4.\medskip

If the labels are ignored, on Pictures  1a) and 1b) one gets the
degree 4 relations due to Sternberg \cite{stern}. The relation on
Picture 1c) is of degree 4, but due to this relation we have to go
twice against the orientation (this relation was indicated by
\cite{stem}). The relations 1b) and 1c) imply the relation on
Picture 3, the degree $5$ relation due to Sternberg \cite{stern}
(if the labels are ignored). One can easy check that the degree 5
relation and the relation 1b) imply the relation 1c). (The dual
degree 5 relation is depicted at Picture 4.)
\medskip

\unitlength=0.60mm \special{em:linewidth 0.4pt}
\linethickness{0.4pt}
\begin{picture}(112.00,65.00)
\put(80.00,10.00){\color{green}\vector(1,0){15.00}}
\put(95.00,10.00){\color{red}\vector(0,1){15.00}}
\put(95.00,25.00){\color{red}\vector(0,1){15.00}}
\put(95.00,40.00){\color{red}\vector(0,1){15.00}}
\put(95.00,55.00){\color{green}\vector(1,0){15.00}}
\put(80.00,10.00){\color{red}\vector(0,1){17.00}}
\put(80.00,27.00){\color{green}\vector(1,0){17.00}}
\put(97.00,27.00){\color{red}\vector(0,1){13.00}}
\put(97.00,40.00){\color{green}\vector(1,0){13.00}}
\put(110.00,40.00){\color{red}\vector(0,1){15.00}}
\put(80.00,10.00){\circle*{2.00}}
\put(111.00,56.00){\circle*{2.00}}
\put(78.00,18.00){\makebox(0,0)[cc]{$1$}}
\put(87.00,7.00){\makebox(0,0)[cc]{$1$}}
\put(97.00,18.00){\makebox(0,0)[cc]{$0$}}
\put(87.00,25.00){\makebox(0,0)[cc]{$\frac 12$}}
\put(93.00,34.00){\makebox(0,0)[cc]{$0$}}
\put(99.00,34.00){\makebox(0,0)[cc]{$0$}}
\put(104.00,38.00){\makebox(0,0)[cc]{$1$}}
\put(112.00,47.00){\makebox(0,0)[cc]{$0$}}
\put(93.00,47.00){\makebox(0,0)[cc]{$1$}}
\put(103.00,57.00){\makebox(0,0)[cc]{$\frac 12$}}
\put(10.00,10.00){\color{green}\vector(1,0){15.00}}
\put(25.00,10.00){\color{red}\vector(0,1){15.00}}
\put(25.00,25.00){\color{red}\vector(0,1){15.00}}
\put(10.00,10.00){\color{red}\vector(0,1){17.00}}
\put(10.00,27.00){\color{green}\vector(1,0){17.00}}
\put(10.00,10.00){\circle*{2.00}}
\put(8.00,18.00){\makebox(0,0)[cc]{$1$}}
\put(17.00,7.00){\makebox(0,0)[cc]{$1$}}
\put(27.00,18.00){\makebox(0,0)[cc]{$0$}}
\put(23.00,32.00){\makebox(0,0)[cc]{$0$}}
\put(17.00,25.00){\makebox(0,0)[cc]{$\frac 12$}}
\put(52.00,30.00){\makebox(0,0)[cc]{$\Rightarrow$}}
\end{picture}
\begin{center}
Picture 3.
\end{center}

\unitlength=0.60mm \special{em:linewidth 0.4pt}
\linethickness{0.4pt}
\begin{picture}(112.00,55.00) \put(80.00,10.00){\color{green}\vector(1,0){15.00}}
\put(95.00,10.00){\color{red}\vector(0,1){15.00}}
\put(95.00,25.00){\color{red}\vector(0,1){15.00}}
\put(95.00,40.00){\color{red}\vector(0,1){15.00}}
\put(95.00,55.00){\color{green}\vector(1,0){15.00}}
\put(80.00,10.00){\color{red}\vector(0,1){17.00}}
\put(80.00,27.00){\color{green}\vector(1,0){17.00}}
\put(97.00,27.00){\color{red}\vector(0,1){13.00}}
\put(97.00,40.00){\color{green}\vector(1,0){13.00}}
\put(110.00,40.00){\color{green}\vector(0,1){15.00}}
\put(80.00,10.00){\circle*{2.00}}
\put(111.00,56.00){\circle*{2.00}}
\put(78.00,18.00){\makebox(0,0)[cc]{$1$}}
\put(87.00,7.00){\makebox(0,0)[cc]{$\frac 12$}}
\put(97.00,18.00){\makebox(0,0)[cc]{$0$}}
\put(87.00,25.00){\makebox(0,0)[cc]{$0$}}
\put(93.00,34.00){\makebox(0,0)[cc]{$1$}}
\put(99.00,34.00){\makebox(0,0)[cc]{$1$}}
\put(104.00,38.00){\makebox(0,0)[cc]{$\frac 12$}}
\put(112.00,47.00){\makebox(0,0)[cc]{$0$}}
\put(93.00,47.00){\makebox(0,0)[cc]{$1$}}
\put(103.00,57.00){\makebox(0,0)[cc]{$0$}}
\put(52.00,30.00){\makebox(0,0)[cc]{$\Rightarrow$}}
\put(25.00,25.00){\color{red}\vector(0,1){15.00}}
\put(25.00,40.00){\color{red}\vector(0,1){15.00}}
\put(25.00,55.00){\color{green}\vector(1,0){15.00}}
\put(27.00,40.00){\color{green}\vector(1,0){13.00}}
\put(40.00,40.00){\color{red}\vector(0,1){15.00}}
\put(41.00,56.00){\circle*{2.00}}
\put(34.00,38.00){\makebox(0,0)[cc]{$\frac 12$}}
\put(42.00,47.00){\makebox(0,0)[cc]{$0$}}
\put(33.00,57.00){\makebox(0,0)[cc]{$0$}}
\put(23.00,33.00){\makebox(0,0)[cc]{$1$}}
\put(23.00,45.00){\makebox(0,0)[cc]{$1$}}
\end{picture}
\begin{center}
Picture 4.
\end{center}

In that follows it will be of use the following two implications
for $B_2$-type Verma relations depicted on Picture 2a, which
follows from implications on Pictures 2, 3 and 4.\bigskip

\unitlength=0.8mm \special{em:linewidth 0.4pt}
\linethickness{0.4pt}
\begin{picture}(128.00,90.00)
\put(78.00,11.00){\circle*{2.00}}
\put(78.00,11.00){\color{green}\vector(1,0){15.00}}
\put(93.00,11.00){\color{red}\vector(0,1){15.00}}
\put(93.00,26.00){\color{red}\vector(0,1){15.00}}
\put(93.00,41.00){\color{green}\vector(1,0){15.00}}
\put(108.00,41.00){\color{red}\vector(0,1){17.00}}
\put(108.00,58.00){\color{green}\vector(1,0){17.00}}
\put(125.00,58.00){\color{red}\vector(0,1){16.00}}
\put(78.00,11.00){\color{red}\vector(0,1){17.00}}
\put(78.00,28.00){\color{green}\vector(1,0){13.00}}
\put(91.00,28.00){\color{red}\vector(0,1){16.00}}
\put(91.00,44.00){\color{green}\vector(1,0){15.00}}
\put(106.00,44.00){\color{red}\vector(0,1){16.00}}
\put(106.00,60.00){\color{red}\vector(0,1){14.00}}
\put(106.00,74.00){\color{green}\vector(1,0){19.00}}
\put(108.00,41.00){\circle*{0.00}}
\put(91.00,44.00){\circle*{2.00}}
\put(125.00,74.00){\circle*{2.00}}
\put(75.00,19.00){\makebox(0,0)[cc]{$1$}}
\put(84.00,31.00){\makebox(0,0)[cc]{$\frac 12$}}
\put(89.00,36.00){\makebox(0,0)[cc]{$0$}}
\put(98.00,46.00){\makebox(0,0)[cc]{$1$}}
\put(104.00,52.00){\makebox(0,0)[cc]{$0$}}
\put(104.00,67.00){\makebox(0,0)[cc]{$1$}}
\put(115.00,77.00){\makebox(0,0)[cc]{$0$}}
\put(128.00,67.00){\makebox(0,0)[cc]{$0$}}
\put(117.00,55.00){\makebox(0,0)[cc]{$\frac 12$}}
\put(110.00,49.00){\makebox(0,0)[cc]{$1$}}
\put(102.00,37.00){\makebox(0,0)[cc]{$0$}}
\put(95.00,32.00){\makebox(0,0)[cc]{$1$}}
\put(95.00,18.00){\makebox(0,0)[cc]{$0$}}
\put(86.00,7.00){\makebox(0,0)[cc]{$1$}}
\put(20.00,0.00){\circle*{2.00}}
\put(20.00,15.00){\color{green}\vector(1,0){15.00}}
\put(35.00,15.00){\color{red}\vector(0,1){15.00}}
\put(50.00,30.00){\color{red}\vector(0,1){15.00}}
\put(50.00,45.00){\color{red}\vector(0,1){15.00}}
\put(20.00,0.00){\color{red}\vector(0,1){15.00}}
\put(35.00,30.00){\color{green}\vector(1,0){15.00}}
\put(47.00,53.00){\makebox(0,0)[cc]{$1$}}
\put(27.00,18.00){\makebox(0,0)[cc]{$\frac 12$}}
\put(43.00,33.00){\makebox(0,0)[cc]{$1$}}
\put(53.00,37.00){\makebox(0,0)[cc]{$0$}}
\put(37.00,22.00){\makebox(0,0)[cc]{$0$}}
\put(23.00,7.00){\makebox(0,0)[cc]{$1$}}
\put(17.00,00.00){\makebox(0,0)[cc]{$A$}}
\put(17.00,16.00){\makebox(0,0)[cc]{$B$}}
\put(37.00,15.00){\makebox(0,0)[cc]{$C$}}
\put(33.00,34.00){\makebox(0,0)[cc]{$D$}}
\put(52.00,30.00){\makebox(0,0)[cc]{$E$}}
\put(52.00,46.00){\makebox(0,0)[cc]{$F$}}
\put(52.00,61.00){\makebox(0,0)[cc]{$G$}}
\put(30.00,65.00){\color{green}\vector(1,0){15.00}}
\put(45.00,65.00){\color{red}\vector(0,1){15.00}}
\put(15.00,20.00){\color{red}\vector(0,1){15.00}}
\put(15.00,35.00){\color{red}\vector(0,1){15.00}}
\put(15.00,50.00){\color{green}\vector(1,0){15.00}}
\put(30.00,50.00){\color{red}\vector(0,1){15.00}}
\put(46.00,80.00){\circle*{2.00}}
\put(12.00,43.00){\makebox(0,0)[cc]{$1$}}
\put(12.00,27.00){\makebox(0,0)[cc]{$0$}}
\put(22.00,47.00){\makebox(0,0)[cc]{$0$}}
\put(27.00,57.00){\makebox(0,0)[cc]{$1$}}
\put(36.00,83.00){\makebox(0,0)[cc]{$0$}}
\put(49.00,73.00){\makebox(0,0)[cc]{$0$}}
\put(38.00,61.00){\makebox(0,0)[cc]{$\frac 12$}}
\put(37.00,45.00){\makebox(0,0)[cc]{$or$}}
\put(63.00,44.00){\makebox(0,0)[cc]{$\Rightarrow$}}
\put(35.00,30.00){\circle*{2.00}}
\end{picture}

\begin{center} Picture 2a.
\end{center}

Let us prove that the bottom implication, that is the graph
$ABCDEFG$ implies the $B_2$-type Verma relation, and we leave to
the reader to prove the upper implication. Because the edge $CD$
is labeled 0, there is a green edge emanating from $C$ (because of
Axiom K3, this edge labelled 1). This implies that there is a
green edge emanating from $A$. Moreover, this edge labeled 1. In
fact if the label would be 0, then, because of Axiom K4, there
would be a commutative diagram, and hence the label of the edge
$BC$ has to be 0, that is not the case. If we assume that there is
the 5 term relation as on the Picture 3, then the label on the
edge $FG$ has to be 0, that is not the case. Hence the $B_2$-type
Verma relation takes place for this graph.\bigskip

{\bf Definition}. A edge-2-colored graph for which Axioms
(K0)-(K5) hold true is called an {\em R-graph} of $B_2$-type.
\medskip

Our aim is to prove the following

\medskip
{\bf Theorem 1}. {\em A connected  R-graph of $B_2$-type  is a
crystal graph of an integrable irreducible
representation of $U_q(sp(4))$, that is a regular $B_2$-type
crystal, and vise versa.}

\section{Construction of R-graphs of $B_2$-type}

Here we give a direct construction of such graphs.  Specifically,
we give a direct constructions of the 3-dimensional "view from the
sky" on a 4-dimensional R-graph $G$.\medskip

{\bf Definition}. Given an R-graph $G$, we define another labeled
digraph $\hat G=(\hat V, \hat E, (X, Y)\colon\hat V\to\mathbb Z_+^2)$ which we call the {\em view from the sky}
(along red strings) on $G$:
\begin{itemize}
\item the set of vertices, $\hat V$, of $\hat G$ is constituted from red strings of $G$, that is a whole red
string in the graph $G$ becomes a vertex in $\hat G$; to each vertex $\hat v$ is attributed a label being a pair
of numbers $(X(\hat v),Y(\hat v))=(t_1(*),h_1(*))$, where we let $*$ to denote the critical
point on the red string $P_1$ corresponding to $\hat v$;

\item the set of edges, $\hat E$, is formed by the following rule: we join vertices $\hat v$ and $\hat v'$ by an
edge going from $\hat v$ to $\hat v'$ and put on it a label either $1$, or $\frac 12$, or $0$, if there exists a
green edge labeled $1$, or $\frac 12$, or $0$, respectively, which joins a pair of vertices (in $G$) located on
the red strings which correspond to $\hat v$ and $\hat v'$, respectively.\end{itemize}

{\bf Lemma 1}. {\em Let $G$ be an R-graph. 
Then any vertex $\hat v$ of $\hat G$ 
has at most three ingoing edges
and at most three outgoing edges; 
the ingoing edges have different labels and the outgoing edges have different
labels; there are no parallel edges in $\hat G$.} \medskip

For a proof of this Lemma we use the following properties, which
follow from the "local commutative diagrams" illustrated on
Pictures 1 and 2 (and dual to them) and the commutative squares
figured in Axiom K4.

\begin{itemize}
\item[(*)] If a vertex $v$ is located at least two red edges above a
critical vertex on $P_1(v)$ or below the critical vertex on
$P_1(v)$, that is $t_1(v)\ge X+2$ or $t_1(v)\le X-2$, respectively. 
Then a pair of the red and green edges ingoing
in $v$ form a commutative square, and a pair of the red and green
edges outgoing from $v$ form  a commutative square.

If a vertex $v$ is such that $t_1(v)=X+1$ (one edge above the
critical), then the pair of ingoing red and green edges form a
commutative square.

If a vertex $v$ is such that $t_1(v)=X-1$ (one edge below the
critical), then the pair of outgoing red and green edges form a
commutative square.
\end{itemize}


\begin{itemize}
\item[(**)]
For each red string $P$, there can be at most two green strings
having critical edge $e$ (i.e., labeled $\frac 12$) such that $e$
enters or leaves a vertex in $P$, and if this is the case, then
the picture is as follows
  \end{itemize}

\unitlength=0.6mm \special{em:linewidth 0.4pt}
\linethickness{0.4pt}
\begin{picture}(95.00,70.00) 
\put(80.00,25.00){\color{red}\vector(0,1){15.00}}
\put(80.00,40.00){\color{red}\vector(0,1){15.00}}
\put(80.00,55.00){\color{red}\vector(0,1){15.00}}
\put(80.00,40.00){\circle*{2.00}}
\put(40.00,40.00){\circle*{2.00}}
\put(80.00,55.00){\color{green}\vector(1,0){15.00}}
\put(83.00,31.00){\makebox(0,0)[cc]{$0$}}
\put(75.00,47.00){\makebox(0,0)[cc]{$1$}}
\put(87.00,58.00){\makebox(0,0)[cc]{$\frac 12$}}
\put(40.00,10.00){\color{red}\vector(0,1){15.00}}
\put(40.00,25.00){\color{red}\vector(0,1){15.00}}
\put(40.00,40.00){\color{red}\vector(0,1){15.00}}
\put(25.00,25.00){\color{green}\vector(1,0){15.00}}
\put(33.00,28.00){\makebox(0,0)[cc]{$\frac 12$}}
\put(43.00,31.00){\makebox(0,0)[cc]{$0$}}
\put(35.00,47.00){\makebox(0,0)[cc]{$1$}}
\end{picture}

\begin{center}Picture 5.\end{center}

That is, either $e$ enters $P$ one red edge before the critical
point or leaves $P$ one red edge after this point.\medskip

{\em Proof of Lemma 1}. Due to the properties (*) and (**), in an R-graph, any red string can be joined with at
most six other red strings through ingoing and outgoing green edges. The claim that the labels are different and
that there are no parallel edges also follow from these properties. \hfill Q.E.D.\medskip

For each R-graph $G$, it will be useful to color edges in the
"view from the sky" graph $\hat G$ in two colors, green and blue.
The edges labeled $\frac 12$ remain green. An edge $\hat e$
labeled 0 (or 1) is colored in blue if $|X(\hat v)-X(\hat v')|=1$
holds, and $\hat e$ remains green if either $|X(\hat v)-X(\hat
v')|=2$ or $|X(\hat v)-X(\hat v')|=0$ hold.

On the following pictures we illustrate "views from the sky" along
a `typical' red strings in an R-graph.\bigskip

\unitlength=0.6mm \special{em:linewidth 0.4pt}
\linethickness{0.4pt}
\begin{picture}(151.00,145.00)
\put(10.00,125.00){\color{green}\vector(1,0){15.00}}
\put(27.00,125.00){\color{green}\vector(1,0){13.00}}
\put(14.00,110.00){\color{green}\vector(3,4){10.67}}
\put(27.00,124.00){\color{green}\vector(3,-4){10.67}}
\put(14.00,141.00){\color{green}\vector(3,-4){11.33}}
\put(26.00,126.00){\color{green}\vector(3,4){11.33}}
\put(26.00,125.00){\color{red}\circle*{2.00}}
\put(9.00,125.00){\color{red}\circle*{2.00}}
\put(42.00,125.00){\color{red}\circle*{2.00}}
\put(14.00,110.00){\color{red}\circle*{2.00}}
\put(38.00,110.00){\color{red}\circle*{2.00}}
\put(14.00,141.00){\color{red}\circle*{2.00}}
\put(38.00,142.00){\color{red}\circle*{2.00}}
\put(8.00,128.00){\makebox(0,0)[cc]{$A$}}
\put(15.00,145.00){\makebox(0,0)[cc]{$B$}}
\put(41.00,145.00){\makebox(0,0)[cc]{$C$}}
\put(45.00,125.00){\makebox(0,0)[cc]{$D$}}
\put(42.00,109.00){\makebox(0,0)[cc]{$E$}}
\put(11.00,108.00){\makebox(0,0)[cc]{$F$}}
\put(26.00,120.00){\makebox(0,0)[cc]{$O$}}
\put(18.00,119.00){\makebox(0,0)[cc]{$0$}}
\put(30.00,135.00){\makebox(0,0)[cc]{$0$}}
\put(21.00,135.00){\makebox(0,0)[cc]{$1$}}
\put(34.00,119.00){\makebox(0,0)[cc]{$1$}}
\put(17.00,127.00){\makebox(0,0)[cc]{$\frac 12$}}
\put(33.00,127.00){\makebox(0,0)[cc]{$\frac 12$}}
\put(120.00,10.00){\color{red}\vector(0,1){15.00}}
\put(120.00,25.00){\color{red}\vector(0,1){15.00}}
\put(120.00,40.00){\color{red}\vector(0,1){15.00}}
\put(120.00,55.00){\color{red}\vector(0,1){15.00}}
\put(120.00,70.00){\color{red}\vector(0,1){15.00}}
\put(120.00,85.00){\color{red}\vector(0,1){15.00}}
\put(120.00,100.00){\color{red}\vector(0,1){15.00}}
\put(120.00,121.00){\makebox(0,0)[cc]{$\vdots$}}
\put(120.00,4.00){\makebox(0,0)[cc]{$\vdots$}}
\put(120.00,125.00){\makebox(0,0)[cc]{$O$}}
\put(120.00,70.00){\circle*{2.00}}
\put(120.00,10.00){\color{green}\vector(1,0){15.00}}
\put(122.00,17.00){\makebox(0,0)[cc]{$0$}}
\put(122.00,32.00){\makebox(0,0)[cc]{$0$}}
\put(122.00,48.00){\makebox(0,0)[cc]{$0$}}
\put(122.00,62.00){\makebox(0,0)[cc]{$0$}}
\put(135.00,10.00){\color{red}\vector(0,1){15.00}}
\put(135.00,25.00){\color{red}\vector(0,1){15.00}}
\put(135.00,40.00){\color{red}\vector(0,1){15.00}}
\put(135.00,55.00){\color{red}\vector(0,1){15.00}}
\put(137.00,17.00){\makebox(0,0)[cc]{$0$}}
\put(137.00,32.00){\makebox(0,0)[cc]{$0$}}
\put(137.00,48.00){\makebox(0,0)[cc]{$0$}}
\put(137.00,62.00){\makebox(0,0)[cc]{$0$}}
\put(120.00,85.00){\color{green}\vector(1,0){15.00}}
\put(135.00,70.00){\color{red}\vector(1,2){6.67}}
\put(141.67,83.00){\color{green}\vector(1,0){0.33}}
\put(142.00,83.00){\color{red}\vector(1,2){6.67}}
\put(148.67,96.00){\color{red}\vector(0,1){15.00}}
\put(141.00,74.00){\makebox(0,0)[cc]{$0$}}
\put(147.00,87.00){\makebox(0,0)[cc]{$0$}}
\put(151.00,102.00){\makebox(0,0)[cc]{$1$}}
\put(127.00,7.00){\makebox(0,0)[cc]{$1$}}
\put(120.00,25.00){\color{green}\vector(1,0){15.00}}
\put(127.00,22.00){\makebox(0,0)[cc]{$1$}}
\put(120.00,40.00){\color{green}\vector(1,0){15.00}}
\put(127.00,37.00){\makebox(0,0)[cc]{$1$}}
\put(120.00,55.00){\color{green}\vector(1,0){15.00}}
\put(127.00,52.00){\makebox(0,0)[cc]{$1$}}
\put(120.00,70.00){\color{green}\vector(1,0){15.00}}
\put(127.00,67.00){\makebox(0,0)[cc]{$1$}}
\put(135.00,85.00){\color{red}\vector(1,2){6.00}}
\put(128.00,87.00){\makebox(0,0)[cc]{$\frac 12$}}
\put(130.00,80.00){\makebox(0,0)[cc]{$0$}}
\put(136.00,91.00){\makebox(0,0)[cc]{$0$}}
\put(102.00,55.00){\color{green}\vector(1,0){18.00}}
\put(111.00,58.00){\makebox(0,0)[cc]{$\frac 12$}}
\put(129.00,73.00){\color{red}\vector(1,2){6.00}}
\put(102.00,55.00){\color{red}\vector(1,2){6.00}}
\put(96.00,43.00){\color{red}\vector(1,2){6.00}}
\put(97.00,50.00){\makebox(0,0)[cc]{$1$}}
\put(103.00,62.00){\makebox(0,0)[cc]{$1$}}
\put(141.00,97.00){\color{red}\vector(0,1){16.00}}
\put(96.00,29.00){\color{red}\vector(0,1){14.00}}
\put(93.00,35.00){\makebox(0,0)[cc]{$0$}}
\put(143.00,105.00){\makebox(0,0)[cc]{$1$}}
\put(104.00,70.00){\color{green}\vector(1,0){15.00}}
\put(104.00,70.00){\color{red}\vector(0,1){15.00}}
\put(104.00,85.00){\color{red}\vector(0,1){15.00}}
\put(104.00,100.00){\color{red}\vector(0,1){15.00}}
\put(104.00,85.00){\color{green}\vector(1,0){16.00}}
\put(104.00,100.00){\color{green}\vector(1,0){16.00}}
\put(104.00,115.00){\color{green}\vector(1,0){16.00}}
\put(120.00,115.00){\color{green}\vector(1,0){14.00}}
\put(120.00,100.00){\color{green}\vector(1,0){14.00}}
\put(134.00,100.00){\color{red}\vector(0,1){15.00}}
\put(105.00,10.00){\color{red}\vector(0,1){15.00}}
\put(105.00,25.00){\color{red}\vector(0,1){15.00}}
\put(105.00,10.00){\color{green}\vector(1,0){15.00}}
\put(107.00,17.00){\makebox(0,0)[cc]{$0$}}
\put(107.00,32.00){\makebox(0,0)[cc]{$0$}}
\put(112.00,7.00){\makebox(0,0)[cc]{$1$}}
\put(105.00,25.00){\color{green}\vector(1,0){15.00}}
\put(112.00,22.00){\makebox(0,0)[cc]{$1$}}
\put(105.00,40.00){\color{green}\vector(1,0){15.00}}
\put(112.00,37.00){\makebox(0,0)[cc]{$1$}}
\put(105.00,40.00){\color{red}\vector(1,2){6.00}}
\put(98.00,58.00){\color{red}\vector(1,2){6.00}}
\put(128.00,88.00){\color{red}\vector(1,2){6.00}}
\put(105.00,40.00){\circle*{2.00}}
\put(148.00,96.00){\circle*{2.00}}
\put(141.00,97.00){\circle*{2.00}}
\put(96.00,43.00){\circle*{2.00}}
\put(134.00,100.00){\circle*{2.00}}
\put(92.00,46.00){\color{red}\vector(1,2){6.00}}
\put(92.00,46.00){\circle*{2.00}}
\put(91.00,31.00){\color{red}\vector(0,1){14.00}}
\put(118.00,77.00){\makebox(0,0)[cc]{$1$}}
\put(118.00,91.00){\makebox(0,0)[cc]{$1$}}
\put(118.00,106.00){\makebox(0,0)[cc]{$1$}}
\put(112.00,113.00){\makebox(0,0)[cc]{$0$}}
\put(126.00,113.00){\makebox(0,0)[cc]{$0$}}
\put(92.00,54.00){\makebox(0,0)[cc]{$1$}}
\put(96.00,22.00){\makebox(0,0)[cc]{$A$}}
\put(89.00,26.00){\makebox(0,0)[cc]{$F$}}
\put(104.00,7.00){\makebox(0,0)[cc]{$B$}}
\put(149.00,116.00){\makebox(0,0)[cc]{$E$}}
\put(142.00,118.00){\makebox(0,0)[cc]{$D$}}
\put(134.00,119.00){\makebox(0,0)[cc]{$C$}}
\end{picture}
\begin{center}
Picture 6a.
\end{center}

\unitlength=0.6mm \special{em:linewidth 0.4pt}
\linethickness{0.4pt}
\begin{picture}(151.00,145.00)
\put(10.00,125.00){\color{green}\vector(1,0){15.00}}
\put(27.00,125.00){\color{green}\vector(1,0){13.00}}
\put(14.00,110.00){\color{blue}\vector(3,4){10.67}}
\put(27.00,124.00){\color{blue}\vector(3,-4){10.67}}
\put(14.00,141.00){\color{green}\vector(3,-4){11.33}}
\put(26.00,126.00){\color{green}\vector(3,4){11.33}}
\put(26.00,125.00){\color{red}\circle*{2.00}}
\put(9.00,125.00){\color{red}\circle*{2.00}}
\put(42.00,125.00){\color{red}\circle*{2.00}}
\put(14.00,110.00){\color{red}\circle*{2.00}}
\put(38.00,110.00){\color{red}\circle*{2.00}}
\put(14.00,141.00){\color{red}\circle*{2.00}}
\put(38.00,142.00){\color{red}\circle*{2.00}}
\put(8.00,128.00){\makebox(0,0)[cc]{$A$}}
\put(15.00,145.00){\makebox(0,0)[cc]{$B$}}
\put(41.00,145.00){\makebox(0,0)[cc]{$C$}}
\put(45.00,125.00){\makebox(0,0)[cc]{$D$}}
\put(42.00,109.00){\makebox(0,0)[cc]{$E$}}
\put(11.00,108.00){\makebox(0,0)[cc]{$F$}}
\put(26.00,120.00){\makebox(0,0)[cc]{$O$}}
\put(18.00,119.00){\makebox(0,0)[cc]{$0$}}
\put(30.00,135.00){\makebox(0,0)[cc]{$0$}}
\put(21.00,135.00){\makebox(0,0)[cc]{$1$}}
\put(34.00,119.00){\makebox(0,0)[cc]{$1$}}
\put(17.00,127.00){\makebox(0,0)[cc]{$\frac 12$}}
\put(33.00,127.00){\makebox(0,0)[cc]{$\frac 12$}}
\put(120.00,10.00){\color{red}\vector(0,1){15.00}}
\put(120.00,25.00){\color{red}\vector(0,1){15.00}}
\put(120.00,40.00){\color{red}\vector(0,1){15.00}}
\put(120.00,55.00){\color{red}\vector(0,1){15.00}}
\put(120.00,70.00){\color{red}\vector(0,1){15.00}}
\put(120.00,85.00){\color{red}\vector(0,1){15.00}}
\put(120.00,100.00){\color{red}\vector(0,1){15.00}}
\put(120.00,121.00){\makebox(0,0)[cc]{$\vdots$}}
\put(120.00,4.00){\makebox(0,0)[cc]{$\vdots$}}
\put(120.00,70.00){\circle*{2.00}}
\put(120.00,10.00){\color{green}\vector(1,0){15.00}}
\put(122.00,17.00){\makebox(0,0)[cc]{$0$}}
\put(122.00,32.00){\makebox(0,0)[cc]{$0$}}
\put(122.00,48.00){\makebox(0,0)[cc]{$0$}}
\put(122.00,62.00){\makebox(0,0)[cc]{$0$}}
\put(135.00,10.00){\color{red}\vector(0,1){15.00}}
\put(135.00,25.00){\color{red}\vector(0,1){15.00}}
\put(135.00,40.00){\color{red}\vector(0,1){15.00}}
\put(135.00,55.00){\color{red}\vector(0,1){15.00}}
\put(137.00,17.00){\makebox(0,0)[cc]{$0$}}
\put(137.00,32.00){\makebox(0,0)[cc]{$0$}}
\put(137.00,48.00){\makebox(0,0)[cc]{$0$}}
\put(137.00,62.00){\makebox(0,0)[cc]{$0$}}
\put(120.00,85.00){\color{green}\vector(1,0){15.00}}
\put(135.00,70.00){\color{red}\vector(1,2){6.67}}
\put(141.67,83.00){\color{green}\vector(1,0){0.33}}
\put(142.00,83.00){\color{red}\vector(1,2){6.67}}
\put(148.67,96.00){\color{red}\vector(0,1){15.00}}
\put(141.00,74.00){\makebox(0,0)[cc]{$0$}}
\put(147.00,87.00){\makebox(0,0)[cc]{$1$}}
\put(151.00,102.00){\makebox(0,0)[cc]{$1$}}
\put(127.00,7.00){\makebox(0,0)[cc]{$1$}}
\put(120.00,25.00){\color{green}\vector(1,0){15.00}}
\put(127.00,22.00){\makebox(0,0)[cc]{$1$}}
\put(120.00,40.00){\color{green}\vector(1,0){15.00}}
\put(127.00,37.00){\makebox(0,0)[cc]{$1$}}
\put(120.00,55.00){\color{green}\vector(1,0){15.00}}
\put(127.00,52.00){\makebox(0,0)[cc]{$1$}}
\put(120.00,70.00){\color{green}\vector(1,0){15.00}}
\put(127.00,67.00){\makebox(0,0)[cc]{$1$}}
\put(135.00,85.00){\color{red}\vector(1,2){6.00}}
\put(128.00,87.00){\makebox(0,0)[cc]{$\frac 12$}}
\put(130.00,80.00){\makebox(0,0)[cc]{$0$}}
\put(136.00,91.00){\makebox(0,0)[cc]{$0$}}
\put(102.00,55.00){\color{green}\vector(1,0){18.00}}
\put(111.00,58.00){\makebox(0,0)[cc]{$\frac 12$}}
\put(129.00,73.00){\color{red}\vector(1,2){6.00}}
\put(102.00,55.00){\color{red}\vector(1,2){6.00}}
\put(96.00,43.00){\color{red}\vector(1,2){6.00}}
\put(97.00,50.00){\makebox(0,0)[cc]{$1$}}
\put(103.00,62.00){\makebox(0,0)[cc]{$1$}}
\put(141.00,97.00){\color{red}\vector(0,1){16.00}}
\put(96.00,29.00){\color{red}\vector(0,1){14.00}}
\put(93.00,35.00){\makebox(0,0)[cc]{$0$}}
\put(143.00,105.00){\makebox(0,0)[cc]{$1$}}
\put(104.00,70.00){\color{green}\vector(1,0){15.00}}
\put(104.00,70.00){\color{red}\vector(0,1){15.00}}
\put(104.00,85.00){\color{red}\vector(0,1){15.00}}
\put(104.00,100.00){\color{red}\vector(0,1){15.00}}
\put(104.00,85.00){\color{green}\vector(1,0){16.00}}
\put(104.00,100.00){\color{green}\vector(1,0){16.00}}
\put(104.00,115.00){\color{green}\vector(1,0){16.00}}
\put(120.00,115.00){\color{green}\vector(1,0){14.00}}
\put(120.00,100.00){\color{green}\vector(1,0){14.00}}
\put(134.00,100.00){\color{red}\vector(0,1){15.00}}
\put(105.00,10.00){\color{red}\vector(0,1){15.00}}
\put(105.00,25.00){\color{red}\vector(0,1){15.00}}
\put(105.00,10.00){\color{green}\vector(1,0){15.00}}
\put(107.00,17.00){\makebox(0,0)[cc]{$0$}}
\put(107.00,32.00){\makebox(0,0)[cc]{$0$}}
\put(112.00,7.00){\makebox(0,0)[cc]{$1$}}
\put(105.00,25.00){\color{green}\vector(1,0){15.00}}
\put(112.00,22.00){\makebox(0,0)[cc]{$1$}}
\put(105.00,40.00){\color{green}\vector(1,0){15.00}}
\put(112.00,37.00){\makebox(0,0)[cc]{$1$}}
\put(105.00,40.00){\color{red}\vector(1,2){6.00}}
\put(98.00,58.00){\color{red}\vector(1,2){6.00}}
\put(128.00,88.00){\color{red}\vector(1,2){6.00}}
\put(105.00,40.00){\circle*{2.00}}
\put(142.00,82.00){\circle*{2.00}}
\put(141.00,97.00){\circle*{2.00}}
\put(96.00,43.00){\circle*{2.00}}
\put(134.00,100.00){\circle*{2.00}}
\put(92.00,46.00){\color{red}\vector(1,2){6.00}}
\put(98.00,58.00){\circle*{2.00}}
\put(91.00,31.00){\color{red}\vector(0,1){14.00}}
\put(118.00,77.00){\makebox(0,0)[cc]{$1$}}
\put(118.00,91.00){\makebox(0,0)[cc]{$1$}}
\put(118.00,106.00){\makebox(0,0)[cc]{$1$}}
\put(112.00,113.00){\makebox(0,0)[cc]{$0$}}
\put(126.00,113.00){\makebox(0,0)[cc]{$0$}}
\put(92.00,54.00){\makebox(0,0)[cc]{$0$}}
\put(96.00,22.00){\makebox(0,0)[cc]{$A$}}
\put(89.00,26.00){\makebox(0,0)[cc]{$F$}}
\put(104.00,7.00){\makebox(0,0)[cc]{$B$}}
\put(149.00,116.00){\makebox(0,0)[cc]{$E$}}
\put(142.00,118.00){\makebox(0,0)[cc]{$D$}}
\put(134.00,119.00){\makebox(0,0)[cc]{$C$}}
\end{picture}
\begin{center}
Picture 6b.
\end{center}

\unitlength=0.6mm \special{em:linewidth 0.4pt}
\linethickness{0.4pt}
\begin{picture}(129.67,107.00)(0,10) \put(30.00,120.00){\circle*{2.00}}
\put(30.00,120.00){\color{blue}\vector(1,1){14.00}}
\put(30.00,120.00){\color{blue}\vector(1,-1){14.00}}
\put(45.00,105.00){\circle*{2.00}}
\put(45.00,135.00){\circle*{2.00}}
\put(30.00,114.00){\makebox(0,0)[cc]{$O$}}
\put(49.00,137.00){\makebox(0,0)[cc]{$B$}}
\put(49.00,104.00){\makebox(0,0)[cc]{$C$}}
\put(15.00,105.00){\circle*{2.00}}
\put(15.00,105.00){\color{blue}\vector(1,1){14.00}}
\put(15.00,135.00){\circle*{2.00}}
\put(15.00,135.00){\color{blue}\vector(1,-1){14.00}}
\put(13.00,137.00){\makebox(0,0)[cc]{$A$}}
\put(12.00,103.00){\makebox(0,0)[cc]{$D$}}
\put(100.00,20.00){\color{red}\vector(0,1){15.00}}
\put(100.00,35.00){\color{red}\vector(0,1){15.00}}
\put(100.00,50.00){\color{red}\vector(0,1){15.00}}
\put(100.00,65.00){\color{red}\vector(0,1){15.00}}
\put(100.00,20.00){\color{green}\vector(1,0){15.00}}
\put(100.00,35.00){\color{green}\vector(1,0){15.00}}
\put(100.00,50.00){\color{green}\vector(1,0){15.00}}
\put(115.00,20.00){\color{red}\vector(0,1){15.00}}
\put(115.00,35.00){\color{red}\vector(0,1){15.00}}
\put(115.00,50.00){\color{red}\vector(1,2){7.67}}
\put(122.67,65.00){\color{red}\vector(1,2){7.00}}
\put(100.00,65.00){\color{green}\vector(1,0){15.00}}
\put(115.00,65.00){\color{red}\vector(0,1){15.00}}
\put(100.00,80.00){\color{green}\vector(1,0){15.00}}
\put(107.67,51.00){\color{red}\vector(1,2){7.00}}
\put(85.00,50.00){\color{red}\vector(0,1){15.00}}
\put(85.00,65.00){\color{red}\vector(0,1){15.00}}
\put(85.00,50.00){\color{green}\vector(1,0){15.00}}
\put(85.00,65.00){\color{green}\vector(1,0){15.00}}
\put(85.00,80.00){\color{green}\vector(1,0){15.00}}
\put(70.00,21.00){\color{red}\vector(1,2){7.67}}
\put(77.67,36.00){\color{red}\vector(1,2){7.00}}
\put(85.00,20.00){\color{red}\vector(0,1){15.00}}
\put(85.00,20.00){\color{green}\vector(1,0){15.00}}
\put(85.00,35.00){\color{green}\vector(1,0){15.00}}
\put(85.67,35.00){\color{red}\vector(1,2){7.00}}
\put(100.00,50.00){\circle*{2.00}}
\put(78.00,36.00){\circle*{2.00}}
\put(85.00,35.00){\circle*{2.00}}
\put(123.00,65.00){\circle*{2.00}}
\put(115.00,65.00){\circle*{2.00}}
\put(94.00,83.00){\makebox(0,0)[cc]{$0$}}
\put(107.00,83.00){\makebox(0,0)[cc]{$0$}}
\put(102.00,72.00){\makebox(0,0)[cc]{$1$}}
\put(102.00,58.00){\makebox(0,0)[cc]{$1$}}
\put(106.00,47.00){\makebox(0,0)[cc]{$1$}}
\put(106.00,15.00){\makebox(0,0)[cc]{$1$}}
\put(93.00,16.00){\makebox(0,0)[cc]{$1$}}
\put(87.00,25.00){\makebox(0,0)[cc]{$0$}}
\put(92.00,32.00){\makebox(0,0)[cc]{$1$}}
\put(92.00,42.00){\makebox(0,0)[cc]{$1$}}
\put(121.00,54.00){\makebox(0,0)[cc]{$0$}}
\put(128.00,70.00){\makebox(0,0)[cc]{$1$}}
\put(70.00,31.00){\makebox(0,0)[cc]{$0$}}
\put(78.00,44.00){\makebox(0,0)[cc]{$1$}}
\put(81.00,58.00){\makebox(0,0)[cc]{$1$}}
\put(81.00,71.00){\makebox(0,0)[cc]{$1$}}
\put(83.00,85.00){\makebox(0,0)[cc]{$D$}}
\put(84.00,13.00){\makebox(0,0)[cc]{$A$}}
\put(116.00,86.00){\makebox(0,0)[cc]{$B$}}
\put(116.00,13.00){\makebox(0,0)[cc]{$C$}}
\end{picture}

\begin{center}
Picture 6c.
\end{center}

\subsection{The "view from the sky" on relations in Axiom K5}
On the next pictures we illustrate the ``view from the sky''
on the relations in R-graphs of $B_2$-type.
\bigskip

\unitlength=0.8mm \special{em:linewidth 0.4pt}
\linethickness{0.4pt}
\begin{picture}(120.00,119.00)(0,20) \put(20.00,135.00){\color{blue}\vector(1,-1){15.00}}
\put(20.00,105.00){\color{blue}\vector(1,1){15.00}}
\put(30.00,130.00){\makebox(0,0)[cc]{$1$}}
\put(30.00,110.00){\makebox(0,0)[cc]{$0$}}
\put(20.00,139.00){\makebox(0,0)[cc]{$(a-1,b-1)$}}
\put(43.00,120.00){\makebox(0,0)[cc]{$(a,b)$}}
\put(20.00,98.00){\makebox(0,0)[cc]{$(a+1,b+1)$}}
\put(20.00,135.00){\circle*{2.00}}
\put(20.00,105.00){\circle*{2.00}}
\put(36.00,120.00){\circle*{2.00}}
\put(28.00,79.00){\makebox(0,0)[cc]{or}}
\put(100.00,95.00){\color{blue}\vector(1,-1){15.00}}
\put(100.00,65.00){\color{blue}\vector(1,1){15.00}}
\put(110.00,90.00){\makebox(0,0)[cc]{$1$}}
\put(110.00,70.00){\makebox(0,0)[cc]{$0$}}
\put(100.00,99.00){\makebox(0,0)[cc]{$(a-1,b-1)$}}
\put(123.00,80.00){\makebox(0,0)[cc]{$(a,b)$}}
\put(100.00,58.00){\makebox(0,0)[cc]{$(a+1,b+1)$}}
\put(100.00,95.00){\circle*{2.00}}
\put(100.00,65.00){\circle*{2.00}}
\put(116.00,80.00){\circle*{2.00}}
\put(85.00,80.00){\color{blue}\vector(1,-1){15.00}}
\put(95.00,75.00){\makebox(0,0)[cc]{$1$}}
\put(85.00,80.00){\circle*{2.00}}
\put(85.00,80.00){\color{blue}\vector(1,1){15.00}}
\put(95.00,85.00){\makebox(0,0)[cc]{$0$}}
\put(79.00,81.00){\makebox(0,0)[cc]{$(a,b)$}}
\put(25.00,64.00){\makebox(0,0)[cc]{$(a-1,b-1)$}}
\put(25.00,23.00){\makebox(0,0)[cc]{$(a+1,b+1)$}}
\put(25.00,60.00){\circle*{2.00}}
\put(25.00,30.00){\circle*{2.00}}
\put(10.00,45.00){\color{blue}\vector(1,-1){15.00}}
\put(20.00,40.00){\makebox(0,0)[cc]{$1$}}
\put(10.00,45.00){\circle*{2.00}}
\put(10.00,45.00){\color{blue}\vector(1,1){15.00}}
\put(20.00,50.00){\makebox(0,0)[cc]{$0$}}
\put(4.00,47.00){\makebox(0,0)[cc]{$(a,b)$}}
\put(56.00,80.00){\makebox(0,0)[cc]{$\Rightarrow$}}
\end{picture}
\begin{center}
Picture 7a
\end{center}

On Picture 7a we illustrated the relation on Picture 1a) and the
dual relation. Let u note that from Axiom K5 also follows the
following relation: two consequent blue edges labeled 0 and 1,
respectively, the roof of the right hand side, also imply the
right-hand side relation on Picture 7a.\medskip

On the next picture we present the "view from the sky" on the
relation on Picture 1b) and its dual (the dual is depicted at the
north-west corner of Picture 7b).

\unitlength=0.8mm \special{em:linewidth 0.4pt}
\linethickness{0.4pt}
\begin{picture}(126.00,131.00)(0,10) \put(20.00,110.00){\circle*{2.00}}
\put(20.00,110.00){\color{green}\vector(1,1){15.00}}
\put(20.00,110.00){\color{green}\vector(3,-1){15.00}}
\put(35.00,126.00){\circle*{2.00}}
\put(35.00,104.00){\circle*{2.00}}
\put(25.00,118.00){\makebox(0,0)[cc]{$0$}}
\put(25.00,105.00){\makebox(0,0)[cc]{$\frac 12$}}
\put(13.00,110.00){\makebox(0,0)[cc]{$(a,b)$}}
\put(35.00,131.00){\makebox(0,0)[cc]{$(a,b-2)$}}
\put(35.00,97.00){\makebox(0,0)[cc]{$(a+1,b-1)$}}
\put(91.00,74.00){\circle*{2.00}}
\put(91.00,74.00){\color{green}\vector(1,1){15.00}}
\put(106.00,89.00){\color{green}\vector(1,-1){15.00}}
\put(91.00,74.00){\color{green}\vector(3,-1){15.00}}
\put(106.00,90.00){\circle*{2.00}}
\put(106.00,68.00){\circle*{2.00}}
\put(122.00,74.00){\circle*{2.00}}
\put(96.00,82.00){\makebox(0,0)[cc]{$0$}}
\put(116.00,82.00){\makebox(0,0)[cc]{$1$}}
\put(96.00,69.00){\makebox(0,0)[cc]{$\frac 12$}}
\put(84.00,74.00){\makebox(0,0)[cc]{$(a,b)$}}
\put(106.00,95.00){\makebox(0,0)[cc]{$(a,b-2)$}}
\put(136.00,73.00){\makebox(0,0)[cc]{$(a+2,b-2)$}}
\put(106.00,61.00){\makebox(0,0)[cc]{$(a+1,b-1)$}}
\put(106.00,69.00){\color{green}\vector(3,1){15.00}}
\put(115.00,68.00){\makebox(0,0)[cc]{$\frac 12$}}
\put(35.00,49.00){\color{green}\vector(1,-1){15.00}}
\put(35.00,50.00){\circle*{2.00}}
\put(35.00,28.00){\circle*{2.00}}
\put(51.00,34.00){\circle*{2.00}}
\put(45.00,42.00){\makebox(0,0)[cc]{$1$}}
\put(35.00,55.00){\makebox(0,0)[cc]{$(a,b-2)$}}
\put(66.00,33.00){\makebox(0,0)[cc]{$(a+2,b-2)$}}
\put(35.00,21.00){\makebox(0,0)[cc]{$(a+1,b-1)$}}
\put(35.00,29.00){\color{green}\vector(3,1){15.00}}
\put(44.00,28.00){\makebox(0,0)[cc]{$\frac 12$}}
\put(31.00,76.00){\makebox(0,0)[cc]{or}}
\put(65.00,74.00){\makebox(0,0)[cc]{$\Rightarrow$}}
\end{picture}
\begin{center}
Picture 7b
\end{center}

On the next picture we illustrated the "view from the sky" on the
relation at Picture 1c) and its dual.\medskip

\unitlength=1.0mm \special{em:linewidth 0.4pt}
\linethickness{0.4pt}
\begin{picture}(108.00,95.00)(0,35) \put(20.00,130.00){\color{green}\vector(1,0){15.00}}
\put(20.00,130.00){\circle*{2.00}}
\put(20.00,130.00){\color{green}\vector(2,-3){14.67}}
\put(20.00,60.00){\color{green}\vector(1,0){15.00}}
\put(36.00,60.00){\circle*{2.00}}
\put(20.00,37.00){\color{green}\vector(2,3){15.33}}
\put(28.00,133.00){\makebox(0,0)[cc]{$\frac 12$}}
\put(25.00,117.00){\makebox(0,0)[cc]{$1$}}
\put(15.00,132.00){\makebox(0,0)[cc]{\small $(a,b)$}}
\put(46.00,132.00){\makebox(0,0)[cc]{\small $(a+1,b-1)$}}
\put(40.00,105.00){\makebox(0,0)[cc]{\small $(a+2,b)$}}
\put(28.00,83.00){\makebox(0,0)[cc]{$or$}}
\put(16.00,61.00){\makebox(0,0)[cc]{\small $(a,b)$}}
\put(47.00,61.00){\makebox(0,0)[cc]{\small $(a+1,b-1)$}}
\put(18.00,33.00){\makebox(0,0)[cc]{\small $(a+1,b+1)$}}
\put(58.00,83.00){\makebox(0,0)[cc]{$\Rightarrow$}}
\put(88.00,92.00){\color{green}\vector(1,0){15.00}}
\put(88.00,92.00){\circle*{2.00}}
\put(88.00,92.00){\color{green}\vector(2,-3){15.30}}
\put(96.00,95.00){\makebox(0,0)[cc]{$\frac 12$}}
\put(93.00,79.00){\makebox(0,0)[cc]{$1$}}
\put(104.00,92.00){\circle*{2.00}}
\put(83.00,92.00){\makebox(0,0)[cc]{\small $(a,b)$}}
\put(116.00,92.00){\makebox(0,0)[cc]{\small $(a+1,b-1)$}}
\put(112.00,69.00){\makebox(0,0)[cc]{\small $(a+2,b)$}}
\put(88.00,69.00){\color{green}\vector(2,3){15.33}}
\put(88.00,69.00){\color{green}\vector(1,0){15.00}}
\put(26.00,63.00){\makebox(0,0)[cc]{$\frac 12$}}
\put(30.00,47.00){\makebox(0,0)[cc]{$0$}}
\put(100.00,83.00){\makebox(0,0)[cc]{$0$}}
\put(80.00,65.00){\makebox(0,0)[cc]{\small $(a+1,b+1)$}}
\put(96.00,72.00){\makebox(0,0)[cc]{$\frac 12$}}
\end{picture}
\begin{center}
Picture 7c
\end{center}

And finally, the view from the sky on Pictures 2 and 2a is
presented as follows: each of the left-hand side combinations of
the green and blue edges implies  the Verma relation, that takes
the form of commuting "views from the sky" on the crystals of the
fundamental integrable modules of $U_q(sp(4))$.\bigskip

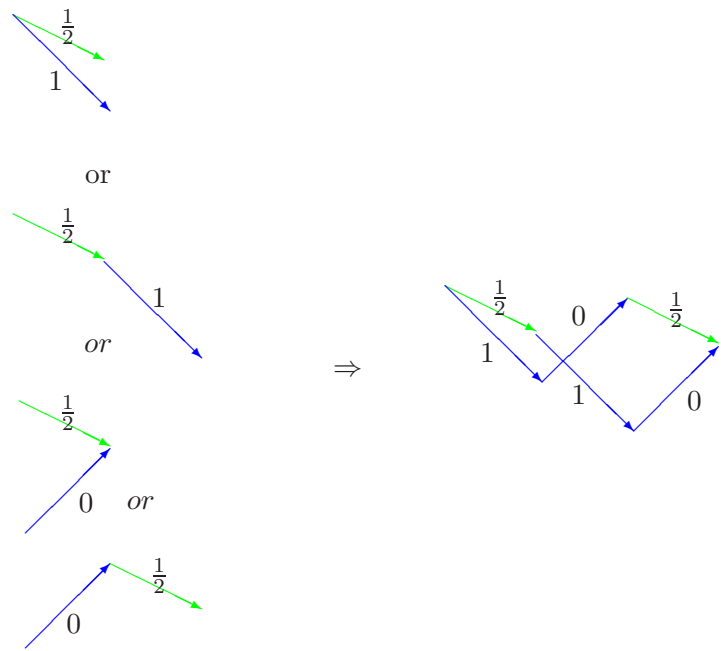
\unitlength=0.80mm \special{em:linewidth 0.4pt}
\linethickness{0.4pt}
\begin{picture}(136.00,100.00)(0,40) \put(91.00,95.00){\color{green}\vector(2,-1){15.00}}
\put(91.00,95.00){\color{blue}\vector(1,-1){16.00}}
\put(106.00,87.00){\color{blue}\vector(1,-1){16.00}}
\put(107.00,79.00){\color{blue}\vector(1,1){14.00}}
\put(122.00,71.00){\color{blue}\vector(1,1){14.00}}
\put(121.00,93.00){\color{green}\vector(2,-1){15.00}}
\put(100.00,93.00){\makebox(0,0)[cc]{$\frac 12$}}
\put(129.00,91.00){\makebox(0,0)[cc]{$\frac 12$}}
\put(98.00,84.00){\makebox(0,0)[cc]{$1$}}
\put(132.00,76.00){\makebox(0,0)[cc]{$0$}}
\put(113.00,77.00){\makebox(0,0)[cc]{$1$}}
\put(113.00,90.00){\makebox(0,0)[cc]{$0$}}
\put(20.00,140.00){\color{green}\vector(2,-1){15.00}}
\put(20.00,140.00){\color{blue}\vector(1,-1){16.00}}
\put(29.00,138.00){\makebox(0,0)[cc]{$\frac 12$}}
\put(27.00,129.00){\makebox(0,0)[cc]{$1$}}
\put(20.00,107.00){\color{green}\vector(2,-1){15.00}}
\put(35.00,99.00){\color{blue}\vector(1,-1){16.00}}
\put(29.00,105.00){\makebox(0,0)[cc]{$\frac 12$}}
\put(22.00,54.00){\color{blue}\vector(1,1){14.00}}
\put(21.00,76.00){\color{green}\vector(2,-1){15.00}}
\put(29.00,74.00){\makebox(0,0)[cc]{$\frac 12$}}
\put(32.00,59.00){\makebox(0,0)[cc]{$0$}}
\put(22.00,35.00){\color{blue}\vector(1,1){14.00}}
\put(36.00,49.00){\color{green}\vector(2,-1){15.00}}
\put(44.00,47.00){\makebox(0,0)[cc]{$\frac 12$}}
\put(44.00,93.00){\makebox(0,0)[cc]{$1$}}
\put(30.00,39.00){\makebox(0,0)[cc]{$0$}}
\put(75.00,81.00){\makebox(0,0)[cc]{$\Rightarrow$}}
\put(34.00,113.00){\makebox(0,0)[cc]{or}}
\put(34.00,85.00){\makebox(0,0)[cc]{$or$}}
\put(41.00,59.00){\makebox(0,0)[cc]{$or$}}
\end{picture}

\begin{center} Picture 7d.
\end{center}

\subsection{Graphs $K(H,A)$}

The above-mentioned graph $K(H,0)$  is an R-graph  of $B_2$-type
which has a red string of the form $(0,H)$ (with $H$ red edges
labeled $1$ and no edges labeled $0$) outgoing from the minimal
vertex and no green edge outgoing from this vertex. These graphs
have the ``view from the sky'' (depending on the parity of $H$)
for $H=1,2,3,4$ as depicted below. One can realize how to
construct such a view for an arbitrary $H$: the building blocks
are the graph $\hat K(2,0)$ (the view from the sky on the
right-hand side relation from Picture 1b (see Picture 7b) and two
rhombuses from the graph $\hat K(3,0)$, these rhombuses are view
from the sky the relations from Picture 3 and 4 (see Pictures 7c).
The rule for travelling along the red edges follows from the
commutativity Axiom K3. We define the set of {\em distinguished
vertices} in $\hat K(H,0)$ to consist of the vertices at the
ground floor, or `zero etage', in the corresponding triangles and
trapezoids; by an analogy with~\cite{a2} we call this set the {\em
diagonal}. Note that the red paths at the diagonal are of the form
$(l,H-l)$, $l=0,\ldots, H$, and the red paths at the etage $k$ are
of the form $(l,H-2k-l)$, $l=0,\ldots, H-2k$. One can see that
$\hat K(H,0)$ has no blue edges. \medskip

\unitlength=0.80mm \special{em:linewidth 0.4pt}
\linethickness{0.4pt}
\begin{picture}(144.00,110.00)
\put(20.00,100.00){\color{green}\vector(1,0){15.00}}
\put(20.00,100.00){\color{red}\circle*{2.00}}
\put(37.00,100.00){\color{red}\circle*{2.00}}
\put(28.00,103.00){\makebox(0,0)[cc]{$\frac 12$}}
\put(18.00,103.00){\makebox(0,0)[cc]{$(0,1)$}}
\put(38.00,103.00){\makebox(0,0)[cc]{$(1,0)$}}
\put(28.00,90.00){\makebox(0,0)[cc]{$\hat K(1,0)$}}
\put(82.00,73.00){\color{green}\vector(1,0){15.00}}
\put(82.00,73.00){\color{red}\circle*{2.00}}
\put(99.00,73.00){\color{red}\circle*{2.00}}
\put(90.00,76.00){\makebox(0,0)[cc]{$\frac 12$}}
\put(77.00,76.00){\makebox(0,0)[cc]{$(0,2)$}}
\put(99.00,76.00){\makebox(0,0)[cc]{$(1,1)$}}
\put(120.00,76.00){\makebox(0,0)[cc]{$(2,0)$}}
\put(100.00,96.00){\makebox(0,0)[cc]{$(0,0)$}}
\put(99.00,73.00){\color{green}\vector(1,0){15.00}}
\put(116.00,73.00){\color{red}\circle*{2.00}}
\put(107.00,76.00){\makebox(0,0)[cc]{$\frac 12$}}
\put(82.00,73.00){\color{green}\vector(1,1){17.00}}
\put(99.00,90.00){\color{green}\vector(1,-1){16.00}}
\put(99.00,91.00){\color{red}\circle*{2.00}}
\put(89.00,84.00){\makebox(0,0)[cc]{$0$}}
\put(108.00,85.00){\makebox(0,0)[cc]{$1$}}
\put(99.00,67.00){\makebox(0,0)[cc]{$\hat K(2,0)$}}
\put(17.00,12.00){\color{green}\vector(1,0){15.00}}
\put(17.00,12.00){\color{red}\circle*{2.00}}
\put(34.00,12.00){\color{red}\circle*{2.00}}
\put(25.00,15.00){\makebox(0,0)[cc]{$\frac 12$}}
\put(34.00,12.00){\color{green}\vector(1,0){15.00}}
\put(51.00,12.00){\color{red}\circle*{2.00}}
\put(42.00,15.00){\makebox(0,0)[cc]{$\frac 12$}}
\put(17.00,12.00){\color{green}\vector(1,1){17.00}}
\put(34.00,29.00){\color{green}\vector(1,-1){16.00}}
\put(34.00,30.00){\color{red}\circle*{2.00}}
\put(24.00,23.00){\makebox(0,0)[cc]{$0$}}
\put(34.00,30.00){\color{green}\vector(1,0){15.00}}
\put(51.00,30.00){\color{red}\circle*{2.00}}
\put(51.00,12.00){\color{green}\vector(1,0){15.00}}
\put(68.00,12.00){\color{red}\circle*{2.00}}
\put(34.00,12.00){\color{green}\vector(1,1){17.00}}
\put(51.00,29.00){\color{green}\vector(1,-1){16.00}}
\put(59.00,24.00){\makebox(0,0)[cc]{$1$}}
\put(37.00,22.00){\makebox(0,0)[cc]{$1$}}
\put(47.00,22.00){\makebox(0,0)[cc]{$0$}}
\put(42.00,33.00){\makebox(0,0)[cc]{$\frac 12$}}
\put(42.00,4.00){\makebox(0,0)[cc]{$\hat K(3,0)$}}
\put(75.00,12.00){\color{green}\vector(1,0){15.00}}
\put(75.00,12.00){\color{red}\circle*{2.00}}
\put(92.00,12.00){\color{red}\circle*{2.00}}
\put(83.00,15.00){\makebox(0,0)[cc]{$\frac 12$}}
\put(92.00,12.00){\color{green}\vector(1,0){15.00}}
\put(109.00,12.00){\color{red}\circle*{2.00}}
\put(100.00,15.00){\makebox(0,0)[cc]{$\frac 12$}}
\put(75.00,12.00){\color{green}\vector(1,1){17.00}}
\put(92.00,29.00){\color{green}\vector(1,-1){16.00}}
\put(92.00,30.00){\color{red}\circle*{2.00}}
\put(82.00,23.00){\makebox(0,0)[cc]{$0$}}
\put(92.00,30.00){\color{green}\vector(1,0){15.00}}
\put(109.00,30.00){\color{red}\circle*{2.00}}
\put(109.00,12.00){\color{green}\vector(1,0){15.00}}
\put(126.00,12.00){\color{red}\circle*{2.00}}
\put(92.00,12.00){\color{green}\vector(1,1){17.00}}
\put(109.00,29.00){\color{green}\vector(1,-1){16.00}}
\put(115.00,24.00){\makebox(0,0)[cc]{$1$}}
\put(95.00,24.00){\makebox(0,0)[cc]{$1$}}
\put(103.00,24.00){\makebox(0,0)[cc]{$0$}}
\put(100.00,33.00){\makebox(0,0)[cc]{$\frac 12$}}
\put(109.00,30.00){\color{green}\vector(1,0){15.00}}
\put(126.00,30.00){\color{red}\circle*{2.00}}
\put(92.00,30.00){\color{green}\vector(1,1){17.00}}
\put(109.00,47.00){\color{green}\vector(1,-1){16.00}}
\put(109.00,48.00){\color{red}\circle*{2.00}}
\put(115.00,41.00){\makebox(0,0)[cc]{$1$}}
\put(100.00,42.00){\makebox(0,0)[cc]{$0$}}
\put(126.00,12.00){\color{green}\vector(1,0){15.00}}
\put(143.00,12.00){\color{red}\circle*{2.00}}
\put(134.00,15.00){\makebox(0,0)[cc]{$\frac 12$}}
\put(109.00,12.00){\color{green}\vector(1,1){17.00}}
\put(126.00,29.00){\color{green}\vector(1,-1){16.00}}
\put(120.00,23.00){\makebox(0,0)[cc]{$0$}}
\put(135.00,24.00){\makebox(0,0)[cc]{$1$}}
\put(118.00,15.00){\makebox(0,0)[cc]{$\frac 12$}}
\put(116.00,33.00){\makebox(0,0)[cc]{$\frac 12$}}
\put(110.00,4.00){\makebox(0,0)[cc]{$\hat K(4,0)$}}
\end{picture}
\begin{center}
Picture 8.
\end{center}

In its turn, the R-graph of $B_2$-type $K(0,A)$ is an  R-graph
which has no red edge outgoing from the minimal vertex, and the
green string beginning at it consists of $A$ edges labeled $1$.
The graph $\hat K(0,A)$ has no green edges. Therefore the
structure of such a graph is forced by the $A_2$-type Verma
relation presented at Picture 7a. This graph  has the form of a
triangular grid, and we depict two examples with $A=1$ and $A=2$
in Picture 9. The distinguished vertex subset in the graph $\hat
K(0,A)$,{ \em the diagonal}, is defined to be constituted by the
unique vertices of the degenerated red paths from the top etage.
\bigskip

\unitlength=0.6mm \special{em:linewidth 0.4pt}
\linethickness{0.4pt}
\begin{picture}(87.00,120.00)
\put(20.00,110.00){\color{red}\circle*{2.00}}
\put(20.00,110.00){\color{blue}\vector(1,-1){15.00}}
\put(36.00,93.00){\color{red}\circle*{2.00}}
\put(37.00,94.00){\color{blue}\vector(1,1){15.00}}
\put(53.00,110.00){\color{red}\circle*{2.00}}
\put(25.00,99.00){\makebox(0,0)[cc]{$1$}}
\put(47.00,99.00){\makebox(0,0)[cc]{$0$}}
\put(20.00,114.00){\makebox(0,0)[cc]{$(0,0)$}}
\put(53.00,114.00){\makebox(0,0)[cc]{$(0,0)$}}
\put(36.00,87.00){\makebox(0,0)[cc]{$(1,1)$}}
\put(36.00,78.00){\makebox(0,0)[cc]{$\hat K(0,1)$}}
\put(20.00,58.00){\color{red}\circle*{2.00}}
\put(20.00,58.00){\color{blue}\vector(1,-1){15.00}}
\put(36.00,41.00){\color{red}\circle*{2.00}}
\put(37.00,42.00){\color{blue}\vector(1,1){15.00}}
\put(53.00,58.00){\color{red}\circle*{2.00}}
\put(25.00,47.00){\makebox(0,0)[cc]{$1$}}
\put(47.00,47.00){\makebox(0,0)[cc]{$0$}}
\put(20.00,62.00){\makebox(0,0)[cc]{$(0,0)$}}
\put(53.00,62.00){\makebox(0,0)[cc]{$(0,0)$}}
\put(53.00,58.00){\color{blue}\vector(1,-1){15.00}}
\put(69.00,41.00){\color{red}\circle*{2.00}}
\put(70.00,42.00){\color{blue}\vector(1,1){15.00}}
\put(86.00,58.00){\color{red}\circle*{2.00}}
\put(58.00,47.00){\makebox(0,0)[cc]{$1$}}
\put(80.00,47.00){\makebox(0,0)[cc]{$0$}}
\put(86.00,62.00){\makebox(0,0)[cc]{$(0,0)$}}
\put(36.00,41.00){\color{blue}\vector(1,-1){15.00}}
\put(52.00,24.00){\color{red}\circle*{2.00}}
\put(53.00,25.00){\color{blue}\vector(1,1){15.00}}
\put(41.00,30.00){\makebox(0,0)[cc]{$1$}}
\put(63.00,30.00){\makebox(0,0)[cc]{$0$}}
\put(32.00,39.00){\makebox(0,0)[cc]{$(1,1)$}}
\put(71.00,38.00){\makebox(0,0)[cc]{$(1,1)$}}
\put(52.00,20.00){\makebox(0,0)[cc]{$(2,2)$}}
\put(52.00,10.00){\makebox(0,0)[cc]{$\hat K(0,2)$}}
\end{picture}
\begin{center}
Picture 9.
\end{center}

Now we are going to give a combinatorial characterization of the
"view from the sky" on an R-graph. A significance of this
characterization is due to that any R-graph is determined by its
"view from the sky".\medskip

The following operation introduced in \cite{a2} will be of use.

Consider arbitrary graphs or digraphs or labeled digraphs $G=(V,E)$ and $H=(V',E')$. Let $S$ be a distinguished
subset of vertices of $G$, and $T$ a distinguished subset of vertices of $H$. Take $|T|$ disjoint copies of $G$,
denoted as $G_t$ ($t\in T$), and $|S|$ disjoint copies of $H$, denoted as $H_s$ ($s\in S$). We glue these copies
together in the following way: for each $s\in S$ and each $t\in T$, the vertex $s$ in $G_t$ is identified with
the vertex $t$ in $H_s$. The resulting graph consisting of $|V| |T|+|V'| |S|- |S| |T|$ vertices and $|E| |T|+|E'|
|S|$ edges is denoted by $(G,S)$\sprod$(H,T)$.

In our case the role of $G$ and $H$ is played by 2-colored
digraphs $\hat K(H,0)$ and $\hat K(0,A)$ depending on parameters
$H,A\in\mathbb Z_+$ (it will be clear later that these digraphs
are the "views from the sky" on the crystals of irreducible
representations of $U_q(sp(4))$ with the highest weight
$H\lambda_1$ and $A\lambda_2$, respectively, where $\lambda_1$ and
$\lambda_2$ denote the fundamental weights for the $B_2$
case).\medskip

Specifically, let us consider a labeled digraph  $\hat K(H,A)=\hat K(H,0)$\sprod $\hat K(0,A)$ (the
distinguished subsets in the graphs $\hat K(H,0)$ and $\hat K(0,A)$ are the diagonals as described above). Then 
we set labels at the vertices of $\hat K(H,A)$ by the following rule: if a vertex $\hat v$ belongs to a
copy of $\hat K(H,0)$, then we set $(\tilde X(\hat v), \tilde Y(\hat v))=(X(\hat v), Y(\hat v))$, where $( X(\hat
v), Y(\hat v))$ is the label on $\hat v$ in $\hat K(H,0)$; if $\hat v$ belongs to a copy of $K(0,A)$ being
attached to the $l$-th distinguished point of $K(H,0)$, then $(\tilde X(\hat v), \tilde Y(\hat v))=(X(\hat v)+l,
Y(\hat v)+H-l)$, where $( X(\hat v), Y(\hat v))$ is the label on $\hat v$ in $\hat K(0,A)$.\bigskip

{\bf Remark}. If we regard $\hat K(H,A)$ as edge green-blue
colored graph, and will follow the rule of changing of the labels
attached to vertices along the edges in according to the color of
the edge and the label attached to the edge, then the labels at
the vertices are determined if we set the label $(0,H)$ to the
source vertex of $\hat K(H,A)$.\medskip

{\bf Definition}. We define $K(H,A)$ to be the R-graph for which
the "view from the sky" is of the form $\hat K(H,A)=\hat
K(H,0)$\sprod $\hat K(0,A)$ with the above defined labels at the vertices.\medskip

{\bf Proposition 1}. {\em Any connected R-graph of $B_2$-type takes the form $K(H,A)$ with $H$ and $A$ being the
lengths of the red and green strings, respectively, emanating from the source vertex.} \medskip

{\em Proof}. In the beginning, we place the list of forbidden
configurations in graphs being the views from the sky on R-graphs on the next Picture.\bigskip

\unitlength=0.8mm \special{em:linewidth 0.4pt}
\linethickness{0.4pt}
\begin{picture}(115.00,85.00)(0,40) \put(30.00,110.00){\color{green}\vector(1,0){15.00}}
\put(45.00,110.00){\color{blue}\vector(1,1){11.00}}
\put(37.00,113.00){\makebox(0,0)[cc]{$\frac 12$}}
\put(52.00,113.00){\makebox(0,0)[cc]{$0$}}
\put(90.00,121.00){\color{blue}\vector(1,-1){11.00}}
\put(101.00,110.00){\color{green}\vector(1,0){14.00}}
\put(97.00,117.00){\makebox(0,0)[cc]{$1$}}
\put(108.00,112.00){\makebox(0,0)[cc]{$\frac 12$}}
\put(41.00,96.00){\makebox(0,0)[cc]{$a)$}}
\put(106.00,96.00){\makebox(0,0)[cc]{$b)$}}
\put(30.00,58.00){\color{blue}\vector(1,1){12.00}}
\put(30.00,58.00){\color{green}\vector(1,0){15.00}}
\put(101.00,70.00){\color{blue}\vector(1,-1){12.00}}
\put(97.00,58.00){\color{green}\vector(1,0){16.00}}
\put(105.00,54.00){\makebox(0,0)[cc]{$\frac 12$}}
\put(38.00,54.00){\makebox(0,0)[cc]{$\frac 12$}}
\put(34.00,67.00){\makebox(0,0)[cc]{$0$}}
\put(108.00,66.00){\makebox(0,0)[cc]{$1$}}
\put(41.00,40.00){\makebox(0,0)[cc]{$c)$}}
\put(106.00,40.00){\makebox(0,0)[cc]{$d)$}}
\end{picture}
\begin{center}
Picture 10
\end{center}
The cases c) and d) easy follow since in these cases there have to
be two outgoing (in c)) or ingoing (in d)) edges for the vertex in
$G$ from which the corresponding green edge labeled $\frac 12$
emanates or terminates, that is not the case due to Axiom K0. The
case a) follows from the $A_2$-type Verma relation depicted on
Picture 1a. The case b) follows the relation depicted on Picture
1b.

Now let $G$ be an R-graph and let us consider a green edge labeled
$\frac 12$. Let us consider the "view from the sky" on the graph
and consider the corresponding edge. Then, because of the "view
from the sky" on the relations on Pictures 7b) and 7c), we get
that such an edge has to be an edge in the graph of the form $\hat
K(H,0)$, for some $H$, we call it a sail of the {\em upper type}.
Now consider a vertex of the bottom etage of this sail. Then,
because of the above list of the forbidden diagrams (Picture 10),
it can be either an ingoing blue edge with the label 0 in such a
vertex or an outgoing blue edge with the label 1, or the both blue
edges, or none. In the latter case we are done. For other cases,
because of the Verma relation and its "view from the sky" we can
translate such a blue edge across all the bottom vertices of the
upper sail and to translate the sail across such an edge. Now with
the help of the $A_2$-type Verma relation (Picture 7a), we get
that the {\em bottom sail} constituted from the blue edges
attached to a bottom vertex takes the "view from the sky" of the
form $\hat K(0,A)$. The point is than any blue edge attached to a
vertex from the below etage of $\hat K(H,0)$, due to Verma
$B_2$-type relation (its view from the sky) implies that a
configuration of the form $\hat K(0, 1)$ is also attached and it
contains such a blue edge. Then applying $A_2$-type Verma
relation, we get the whole bottom sail. Thus, any connected
R-graph $G$ has the "view from the sky" of the form $\hat G=\hat
K(H,0)$\sprod $\hat K(0,A)$. It is clear that $H$ and $A$ are the
lengths of the red and green strings passing through the minimal
vertex of $G$, and the source vertex has the label $(0,H)$.

Vice versa, given a graph $\hat G=\hat K(H,0)$\sprod $\hat K(0,A)$
we reconstruct $G$ using the Pictures 6a, 6b, 6c.\hfill  Q.E.D.\medskip

The rest of the paper is devoted to a proof that the graph
$K(H,A)$ is the crystal graph of the irreducible representation of
$U_q(sp(4))$ with the highest weight $H\lambda_1+A\lambda_2$. To
prove this, we will use a new model related to $B_2$-type
crystals.

\section{Crossings model for $B_2$-type crystals}

Here we present a model of the {\em free regular crystal}
$K^{\infty}$ of $B_2$-type. That is an edge-2-colored graph with
infinite monochromatic strings passing through each vertex, and we
will show that each $B_2$-type crystal for irreducible
representation can be obtained as an interval of $K^{\infty}$.

The vertices of $K^{\infty}$ are the functions on the vertex-set
the diagram depicted on Picture 11, for which the inequalities
indicated by arrows are imposed, namely:  $f(a)\ge f(b)\ge f(c)$
and $f(x)\ge f(y)\ge f(z)\ge f(w)$. Also it is required that the
following parity conditions hold: $f(y)+f(z)\in 2\mathbb Z$ and
$f(a),\,f(c)\in 2\mathbb Z$, and that at least {\em three} of the above inequalities turn into equalities. The
possible combinations of equalities are given in Picture 12.
(One may consider an equivalent model, with odd $f(a),\,f(c)\in 2\mathbb Z+1$ and
half-integer weights.)

Such functions are called {\em admissible configurations}.\medskip

\unitlength=0.6mm \special{em:linewidth 0.4pt}
\linethickness{0.4pt}
\begin{picture}(119.00,60.00)
\put(30.00,10.00){\circle*{2.00}}
\put(30.00,10.00){\vector(1,2){19.67}}
\put(51.00,50.00){\circle*{2.00}}
\put(51.00,50.00){\vector(1,0){38.00}}
\put(91.00,50.00){\circle*{2.00}}
\put(91.00,49.00){\vector(1,-2){19.00}}
\put(111.00,10.00){\circle*{2.00}}
\put(34.00,50.00){\circle*{2.00}}
\put(34.00,50.00){\vector(1,-1){38.00}}
\put(73.00,10.00){\circle*{2.00}}
\put(74.00,11.00){\vector(1,1){38.00}}
\put(114.00,50.00){\circle*{2.00}}
\put(32.00,7.00){\makebox(0,0)[cc]{$x$}}
\put(51.00,54.00){\makebox(0,0)[cc]{$y$}}
\put(91.00,54.00){\makebox(0,0)[cc]{$z$}}
\put(109.00,7.00){\makebox(0,0)[cc]{$w$}}
\put(30.00,52.00){\makebox(0,0)[cc]{$a$}}
\put(73.00,5.00){\makebox(0,0)[cc]{$b$}}
\put(119.00,50.00){\makebox(0,0)[cc]{$c$}}
\end{picture}

\begin{center} Picture 11.
\end{center}

\unitlength=0.8mm \special{em:linewidth 0.4pt}
\linethickness{0.4pt}
\begin{picture}(134.00,131.00)
\put(20.00,110.00){\circle*{2.00}}
\emline{20.00}{110.00}{1}{30.00}{130.00}{2}
\emline{30.00}{130.00}{3}{50.00}{130.00}{4}
\emline{50.00}{130.00}{5}{60.00}{110.00}{6}
\put(60.00,110.00){\circle*{2.00}}
\put(50.00,130.00){\circle*{2.00}}
\put(30.00,130.00){\circle*{2.00}}
\put(23.00,90.00){\circle*{2.00}}
\emline{23.00}{90.00}{7}{40.00}{69.00}{8}
\emline{40.00}{69.00}{9}{57.00}{90.00}{10}
\put(57.00,90.00){\circle*{2.00}}
\put(40.00,69.00){\circle*{2.00}}
\put(50.00,90.00){\circle*{2.00}}
\emline{50.00}{90.00}{11}{60.00}{69.00}{12}
\put(60.00,69.00){\circle*{2.00}}
\put(95.00,90.00){\circle*{2.00}}
\emline{95.00}{90.00}{13}{112.00}{69.00}{14}
\emline{112.00}{69.00}{15}{129.00}{90.00}{16}
\put(129.00,90.00){\circle*{2.00}}
\put(112.00,69.00){\circle*{2.00}}
\put(103.00,90.00){\circle*{2.00}}
\emline{103.00}{90.00}{17}{92.00}{69.00}{18}
\put(92.00,69.00){\circle*{2.00}}
\put(20.00,20.00){\circle*{2.00}}
\emline{20.00}{20.00}{19}{30.00}{40.00}{20}
\emline{30.00}{40.00}{21}{50.00}{40.00}{22}
\put(50.00,40.00){\circle*{2.00}}
\put(30.00,40.00){\circle*{2.00}}
\emline{40.00}{19.00}{23}{57.00}{40.00}{24}
\put(57.00,40.00){\circle*{2.00}}
\put(40.00,19.00){\circle*{2.00}}
\emline{103.00}{40.00}{25}{123.00}{40.00}{26}
\emline{123.00}{40.00}{27}{133.00}{20.00}{28}
\put(133.00,20.00){\circle*{2.00}}
\put(123.00,40.00){\circle*{2.00}}
\put(103.00,40.00){\circle*{2.00}}
\put(95.00,40.00){\circle*{2.00}}
\emline{95.00}{40.00}{29}{112.00}{19.00}{30}
\put(112.00,19.00){\circle*{2.00}}
\put(40.00,101.00){\makebox(0,0)[cc]{$(a)$}}
\put(44.00,56.00){\makebox(0,0)[cc]{$(b1)$}}
\put(110.00,56.00){\makebox(0,0)[cc]{$(b2)$}}
\put(44.00,9.00){\makebox(0,0)[cc]{$(c1)$}}
\put(110.00,9.00){\makebox(0,0)[cc]{$(c2)$}}
\end{picture}

\begin{center} Picture 12.
\end{center}

Now we define the structure of a crystal on this set of admissible
configurations $f$ by explaining to which admissible configuration
$f$ goes under the action of the operations that we denote by
$F_1$ and $F_2$ as before. That is, the red edges of the crystal
take the form $(f,F_1f)$ and the green edges take the form
$(f,F_2f)$.
\medskip

The operation $F_1$ ``moves'' $f$ to an admissible configuration
$f'$ which can differ from $f$ at exactly one of the points  $x$,
$b$, $w$. Specifically,

if $f(w) < f(z)$, then $F_1$ increases $f$ by one at the vertex
$w$: $f'(w)=f(w)+1$;

if $f(w)=f(z)$ and $f(b)<f(a)$, then $F_1$ increases $f$ by one at
$b$: $f'(b)=f(b)+1$;

and if $f(w)=f(z)$ and $f(b)=f(a)$, then $F_1$ increases $f$ by
one at $x$:  $f'(x)=f(x)+1$.
\medskip

In its turn, $f'=F_2f$ differs from $f$ at one or two points among
$a$, $y$, $z$, $c$, and $F_2$ either increases $f$ by 2 at one
point, or increases $f$ by 1 at $y$ and at $z$. More precisely,
\medskip

if  $f(a)- f(b)<f(a)- f(c)$, then $f'(c)= f(c)+2$;

if $f(a)-f(b)\ge f(b)-f(c)$ and $f(y)+2\le f(x)$, then
$f'(y)=y+2$;

if $f(a)-f(b)\ge f(b)-f(c)$ and $f(y)+1= f(x)$, then
$f'(y)=f(y)+1$ and $f'(z)=f(z)+1$;

if $f(a)-f(b)\ge f(b)-f(c)$, $f(y)= f(x)$, and $f(z)<f(y)$, then
$f'(z)=f(z)+2$;

and if $f(a)-f(b)\ge f(b)-f(c)$, $f(y)= f(x)$, and $f(z)=f(y)$,
then $f'(a)=f(a)+2$.
\medskip

It is easy to see that these operations preserve the
admissibility.
\medskip

{\bf Definition}. An admissible configuration is called a {\em fat
vertex} if $f(a)=f(b)=f(c)\in 2\mathbb Z$ and
$f(x)=f(y)=f(z)=f(w)$.
\medskip

The subset $B(H,A)\subset K^{\infty}$ of configurations which
satisfy the restrictions $f(c)\ge 0$, $f(w)\ge 0$, $f(a)\le A$,
$f(x)\le H$ is called the {\em interval} of weight
$A/2\lambda_2+H\lambda_1$ ($A\in 2\mathbb Z$).

{\bf Remark}. 
Note that the interval `join' the fat vertex $\bf 0$ and the fat vertex $f(a)=f(b)=f(c)=A$ and
$f(x)=f(y)=f(z)=f(w)=H$, and the whole rectangle of the fat vertices $0\le f(a)=f(b)=f(c)\le A$ and $0\le
f(x)=f(y)=f(z)=f(w)\le H$ belong to this interval. Also note that the interval joining a fat vertex
$f(a)=f(b)=f(c)=A'$ and $f(x)=f(y)=f(z)=f(w)=H'$ and the fat vertex $f(a)=f(b)=f(c)=A+A'$ and
$f(x)=f(y)=f(z)=f(w)=H+H'$, $A',\, A\in 2\mathbb Z$, $H,\,H'\in\mathbb Z$, $A$, $H\ge 0$, is isomorphic to the
interval of weight $A/2\lambda_2+H\lambda_1$. \medskip

We define the operations $F_1$ and $F_2$ on the interval $B(H,A)$
as they are defined in $K^{\infty}$ with the following
modification at the final cases: if in the last case of the
definition of $F_1$, one has $f(x)=H$, then we say that $F_1$ is not defined at $f$; and if in the last case of
the definition of $F_2$, one has $f(a)=A$, then we say that $F_2$ is not defined at $f$. Accordingly, the
reverse operation $F_2^{-1}$ does not act if it would result in the value on $c$ below zero, and similarly for
$F_1^{-1}$ and $w$. \medskip

\noindent One can check that this model gives the inclusion

\hspace{2cm} $B(H,A)\supset B(H,0)$\sprod $B(0,A)$,

\noindent where the distinguished subsets in $B(H,0)$ and $B(0,A)$
are constituted by the fat points.
\bigskip

In the next section we will prove, using Littelmann's path model
\cite{litt}, that the interval $B(H,A)\subset K^{\infty}$ is a
crystal graph of the integrable irreducible representation of
$U_q(sp(4))$ with weight $\frac {A}{2}\lambda_2+H\lambda_1$.
\medskip

{\bf Remark}. On the free regular crystal $K^{\infty}$, there is
the following involution which reverse the direction of edges 
$$*:K^{\infty}\to K^{\infty},$$ $*f(a)=-f(c)$, $*f(b)=-f(b)$, $*f(c)=-f(a)$, $*f(x)=-f(w)$,
$*f(y)=-f(z)$, $*f(z)=-f(y)$, $*f(w)=-f(x)$. One can see that $*f$ is an admissible configuration, moreover, if
$f$ is the type (a), then $*f$ is also of type (a); if $f$ is of type ($bi$), $i=1,2$, then $*f$ is of type
($b(3-i)$); and if $f$ is of type ($ci$), $i=1,2$, then $*f$ is of type ($c(3-i)$). This involution sends each
interval $B(H,A)$ to an interval $*B(H,A)$. The interval $*B(H,A)$ is isomorphic to the interval $B(H,A)$
(see Remark above), and the composition $B(H,A)\to *B(H,A)\cong B(H,A)$ is the Kashiwara involution
(\cite{kas-95}) indeed.

\section{Littelmann cones and admissible configurations} 
\subsection{Canonical coordinates for configurations}

It is easy to check that for any configuration  $f\in B(H,A)$
(where $f(c)\ge 0$, $f(w)\ge 0$, $f(a)\le A$ and $f(x)\le H$), we
can reach the sink (or the maximal vertex, that is, the fat vertex
with $A=f(a)=f(b)=f(c)$ and $H=f(x)=f(y)=f(z)=f(w)$) by applying a
word of the form $F_1^nF_2^mF_1^pF_2^q$ or
$F_2^{q'}F_1^{p'}F_2^{m'}F_1^{n'}$ for appropriate $n,m,p,q$, and
reach the source (or the minimal vertex, that is, the fat vertex
with $0=f(a)=f(b)=f(c)$ and $0=f(x)=f(y)=f(z)=f(w)$) by applying a
words of the same form but substituting $F_i$ on $E_i$, $i=1,2$.
These words are called {\em
canonical words} for $f$, and the corresponding $(n,m,p,q)$ and
$(q',p',m',n')$ are called the {\em canonical coordinates} of $f$.

We can explicitly calculate these coordinates for an admissible
configuration $f$. Indeed, there are 5 cases, depending on the type of the
equalities diagram. The answers are as follows:\medskip

(a) Let  $\alpha=f(a)$, $\beta=f(b)$, $\gamma=f(c)$, and
$h=f(x)=f(y)=f(z)=f(w)$. Then the following canonical words move
$f$ to the sink. The first word can be of two types, relative to
$\beta$ and $\frac{\alpha+\gamma}2$. Specifically, if
$f(a)- f(b)\ge f(b)-f(a)$, then the word is
  $$
1^{H-h}\,2^{\frac{A-\gamma}2+H-h}\,1^{H-h+A-\beta}\,
2^{\frac{A-\alpha}2},
  $$
  that is $q=\frac{A-\alpha}2$, $p=H-h+A-\beta$,
  $m=\frac{A-\gamma}2+H-h$, $n=H-h$;\\
and if  $f(a)-f(b)<f(b)-f(c)$, then the word is
  $$
1^{H-h}\,2^{\frac{A+\alpha}2-\beta +H-h}\,1^{H-h+A-\beta}\,
2^{\frac{A-\alpha}2+\beta-\frac{\alpha+\gamma}2},
$$
that is $q=\frac{A-\alpha}2+\beta-\frac{\alpha+\gamma}2$,
$p=H-h+A-\beta$, $m=\frac{A+\alpha}2-\beta +H-h$, $n=H-h$

The second word is
  $$
2^{\frac{A-\alpha}2}\,1^{H-h+A-\alpha}\,2^{\frac{A-\gamma}2+H-h}\,
1^{\alpha-\beta+H-h},
  $$
that is $n'=\alpha-\beta+H-h$, $m'=\frac{A-\gamma}2+H-h$,
$q'=H-h+A-\alpha$, $p'=\frac{A-\alpha}2$.

(b1) Let  $\alpha=f(a)=f(b)=f(c)$, $h=f(x)$, $g=f(y)$,
$e=f(z)=f(w)$. Then the first and second words are (respectively)
   $$
1^{H-h}\,2^{H-h+\frac{A-\alpha}2}\,1^{h-e+A-\alpha+H-h}\,2^{h-\frac{g+e}2+
\frac{A-\alpha}2}=
  $$
  $$
1^{H-h}\,2^{H-h+\frac{A-\alpha}2}\,1^{H-e+A-\alpha}\,2^{h-\frac{g+e}2+
\frac{A-\alpha}2};
  $$

 $$2^{\frac{A-\alpha}2}\,1^{H-e+A-\alpha}\,2^{H-\frac{g+e}2+
\frac{A-\alpha}2}\,1^{H-h}.
$$

(b2) Let  $\alpha=f(a)=f(b)=f(c)$, $h=f(x)=f(y)$, $g=f(z)$, $e
=f(w)$. Then the words are
   $$
1^{H-h}\,2^{H-h+\frac{A-\alpha}2 }\,1^{h-e+A-\alpha+H-h}\,2^{\frac{h-g}2+
\frac{A-\alpha}2 };
  $$

  $$
2^{\frac{A-\alpha}2}\,1^{H-g+A-\alpha}\,2^{H-\frac{g+h}2+\frac{A-\alpha}2}\,
1^{g-e+H-h}.
  $$

(c1) Let  $\alpha=f(a)$, $\beta=f(b)=f(c)$, $h=f(x)=f(y)=f(z)$,
$e=f(w)$. Then the words are
  $$
1^{H-h}\,2^{\frac{A-\beta}2+H-h}\,1^{h-e+A-\beta+H-h}\,2^{\frac{A-\alpha}2};
  $$

  $$
2^{\frac{A-\beta}2-[\frac{\alpha-\beta}2]}\,1^{H-h+A-\alpha}\,
2^{[\frac{\alpha-\beta}2]+H-h+\frac{A-\alpha}2}\,1^{h-e+\alpha-\beta+H-h}.
$$

(c2) Let   $\alpha=f(a)=f(b)$, $\gamma=f(c)$, $h=f(x)$,
$g=f(y)=f(z)=f(w)$. Then the words are
   $$
1^{H-h}\,2^{\frac{A-\alpha}2+H-h}\,
1^{h-g+A-\alpha+H-h}\,2^{\frac{\alpha-\gamma}2+h-g+\frac{A-\alpha}2};
  $$

  $$
2^{\frac{A-\alpha}2}\,1^{H-g+A-\alpha}\,2^{\frac{\alpha-\gamma}2+H-g+
\frac{A-\alpha}2}\,1^{H-h}.
$$

\subsection{The cones of canonical words}

According to~\cite{litt}, the canonical words of $B_2$-type
constituted the cones of the following form:
 $$
C_1=\{(n',m'p',q')\,:\,2m'\ge p'\ge 2q'\}\mbox{ and }
C_2=\{(q,p,m,n)\,:\,p\ge m\ge n\}.
  $$
Let us check that the canonical coordinates of admissible
configurations belong to these cones. We have to check gradually
all 5 cases:

  \begin{itemize}

\item[(a)] $C_1$: \ $2(\frac{A-\gamma}2+H-h)\ge H-h+A-\alpha\ge
2\frac{A-\alpha}2$;
  \item[] $C_2$:\ if $2\beta\le\alpha+\gamma$, then there holds
$H-h+A-\beta\ge H-h+\frac{A-\gamma}2\ge H-h$, since $2\beta\le
A+\gamma$.
If $2\beta>\alpha+\gamma$, then there holds $H-h+A-\beta\ge
H-h+\frac{A+\alpha}2-\beta\ge H-h$.

\item[(b1)] $C_1$: \ $2H-g-e+A-\alpha\ge H-e+A-\alpha\ge
2\frac{A-\alpha}2$;
  \item[]
$C_2$:\ $H-e+A-\alpha\ge H-h+ \frac{A-\alpha}2\ge H-h$.

\item[(b2)] $C_1$: \ $2H-g-h+A-\alpha\ge H-g+A-\alpha\ge
2\frac{A-\alpha}2$;
  \item[]
$C_2$:\ $H-e+A-\alpha\ge H-h+ \frac{A-\alpha}2\ge H-h$.

\item[(c1)] $C_1$: \ $2(H-h)+A-\alpha+2\frac{\alpha-\beta}2\ge
H-h+A-\alpha\ge A-\beta-[\frac{\alpha-\beta}2]$;
  \item[]
$C_2$:\ $H-h+A-\beta+h-e\ge H-h+\frac{A-\beta}2\ge H-h$.

\item[(c2)] $C_1$: \ $2(H-g)+A-\gamma\ge H-g+A-\alpha\ge A-\alpha$;
   \item[]
$C_2$:\ $H-h+A-\beta+h-e\ge H-h+\frac{A-\beta}2\ge H-h$.
   \end{itemize}

Thus, the coordinates of the admissible configurations belong to
the cones. Moreover, one can see that configurations of $B(H,A)$
cover the `intervals' $q\le \frac A2$, $p\le A+H$, $m\le H+\frac
A2$, $n\le H$ and $n'\le H$, $m'\le H+\frac A2$, $p'\le A+H$,
$q'\le \frac A2$, respectively.

In order to get the equivalence of our model and the Littelmann
path model, we have to check that the piece-wise linear
transformations on the canonical coordinates of admissible
configurations are the same as  for the Littlemann cones. That is,
we have to check the validity of the following relations:
$$
q=\max (q',p'-m',m'-n');
$$
$$
p=\max (p', n'+2p'-2m', 2q'+n');
$$
$$
m=\min (m', 2m'-p'+q', n'+q');
$$
$$
n=\min (n', 2m'-p',p'-2q').
$$

Note that since
   $$
n+p=n'+p'\mbox{ and } m+q=m'+q',
  $$
it suffices to check the first two equalities.

The most involved case to be checked is the case of the
configurations with the equalities of type (a), and here are the
calculations for this case (other cases are much easier and we
leave them for the reader).

So, in the case of equalities of type (a):

\medskip
If $\beta\le\frac{\alpha+\gamma}2$, then we have
$q=\frac{A-\alpha}2$, $q'=\frac{A-\alpha}2$,
$p'-m'=A-\alpha-\frac{A-\gamma}2$, and
$m'-n'=\frac{A-\gamma}2-(\alpha-\beta)$. So we obtain $q'\ge
m'-n'$, whence  $q=q'$;

\medskip
If $\beta\le\frac{\alpha+\gamma}2$, then
$q=\frac{A-\gamma}2-(\alpha-\beta)$, $q'=\frac{A-\alpha}2$,
$p'-m'=A-\alpha-\frac{A-\gamma}2$, and
$m'-n'=\frac{A-\gamma}2-(\alpha-\beta)$. So we obtain $q'\le
m'-n'$, whence  $q=m'-n'$.

\medskip
Thus, the relation $q=\max (q',p'-m',m'-n')$ is valid.
\medskip

The equality  $ p=\max (p', n'+2p'-2m', 2q'+n') $ holds,
since in the case under consideration,
we have $p=n'+2q'=\max (p', n'+2p'-2m', 2q'+n')$. \medskip

Thus, we can conclude with the following
\medskip

{\bf Theorem 2}. {\em For any $A\in 2\mathbb Z_+$ and $H\ge 0$,
the interval $B(H,A)$ is the crystal graph of an integrable irreducible representation of the rank 2 algebra of
$B_2$-type and vice versa.}\medskip

{\bf Remark}. We want to stress the following aspect of the
crossing model. From this model one can see that a crystal graph
of an integrable representation of $B_2$-type is located on
a union of 5 polyhedra, and the latter set is not a polyhedron itself. There are two projections of this union
of polyhedra to the Littelmann cones, and these cones have to be related via the specific piece-wise linear
transformation. Thus, the crossing model captures the global non-convex structure of regular crystals,  and  the
language of the Littelmann cones describes the projections of this non-convex structure.

\section{$B(H,A)=K(H,A/2)$}

In view of Theorem 2, Theorem 1 would follow from the
equality $B(H,A)=K(H,A/2)$.

In the beginning, we establish the required bijection for the case
$H=0$ and for the case $A=0$.

\medskip
{\bf Claim}. $B(H,0)=K(H,0)$.
\medskip

\noindent {\bf Proof.} Note that the admissible configurations in
the interval $B(H,0)$ have equalities of types (b1) or  (b2).

Take a vertex $v$ in $K(H,0)$. The corresponding admissible
configuration is either of type (b1) or of type (b2). Consider the
red-colored path passing through this vertex. Then this vertex is
located either below the critical vertex (Case 1) or strictly
above it (Case 2). Now consider the corresponding
vertex for this path in the graph being the ``view from the sky'' for
$K(H,0)$. Let this vertex be located on the $k$-th etage and at
$l$ steps to the right (that is, the corresponding red
path is of the form $(l,H-2k-l)$, $l$ edges with label 0 and
$H-2k-l$ with label 1).

\medskip
\noindent{\em Case 1.} Suppose $l=2t$.  We consider the following route from the
source to $v$: first we traverse $2(k+t)$ red edges (these edges have label 1),
then $k+t$ green edges with label 0 and $t$ green edges with label 1,
and finally, $l'\le l=2t$ red edges with label $0$ until we reach
the vertex $v$.

The configuration corresponding to $v$ is of type (b2), and it is
represented as follows: $f(a)=f(b)=f(c)=0$,  $f(x)=2(k+t)$ (= the
number of red edges with label 1 on the route); $f(y)=2(k+t)$ (=
twice the number of green edges with label 0 on the route);
$f(z)=2t+0$ (= twice the number of green edges with label 1 plus 0
edges with label $\frac 12$); and $f(w)=l'\le 2t$ (= the number of
final red edges with label $0$). Conversely, by a configuration of
(b2)-type with $f(z)\in 2\mathbb Z$, we reconstruct the vertex $v$
(in general case, lying in the {\em sail} $f(a)=f(b)=f(c)=\alpha$)
as follows:
  $$
k=\frac{f(x)-f(z)}2,\quad l=f(z),
  $$
and $f(w)$ determines the position of $v$ in the red path $P_1(v)$
(that is $f(w)$ equals the number of edges to go from the
beginning of the path $P_1(v)$ to $v$).
\medskip

Now let $l=2t+1$. We start from the source and traverse one
red edge with label 1, then one green edge with label $\frac 12$, then
$2(k+t)$ red edges (all but one have label 1), then $(k+t)$ green edges
with label 0 and $t$ green edges with label $1$, and finally,
$l'\le 2t+1$ red edges until we reach $v$.

The corresponding configuration is again of type (b2) and represented as:
$f(a)=f(b)=f(c)=0$; $f(x)=1+2(k+t)-1=2(k+t)$ (= the number of
red edges with the label 1); $f(y)=2(k+t)$ (= the number
of green edges with label 0); $f(z)=2t+1$ (= twice the number
of green edges with label 1 plus the number of green edges
with label $\frac 12$); $f(w)=l'\le l=2t+1$ (= the number of
red edges with label $0$ on the final part of the
route). Conversely, a configuration of (b2)-type with $f(z)\in
2\mathbb Z+1$ reconstructs the vertex $v$ as follows:
  $$
k=\frac{f(x)-f(z)-1}2,\quad  l=f(z),
  $$
and $f(w)$ determines   the position of $v$ in the
red path $P_1(v)$.
\medskip

\noindent{\em Case 2.}
Suppose $v$ is located on $P_1(v)$ above the critical point, and
$l=2t$. Then, starting from the source, we traverse $2(k-1)+t$
red edges with label 1, then green $k+t-1$ edges with label 0 and
$t$ green edges with label 1, then move along the red path
until we reach the tail of the green edge with label 0 that
enters $v$ (on this part of the route we have passed through
$2t$ edges with label 0 and $s\ge 1$ edges with label 1), and finally,
traverse this green edge to reach $v$.

The corresponding configuration is of type (b1) for which:
$f(a)=f(b)=f(c)=0$; $f(x)=2(k+t)-2+(s+2)$ (= the number of red
edges with label 1 (note: on the red path from which we turn to
the $v$ through the green edge, the number of edges with label 1
is greater by two; see the corresponding relation on Picture 3,
see also Picture 10a)); $f(y)=2(k+t-1+1)$ (= the number of green
edges with label 0); $f(z)=f(w)=2t$ (= twice the number of green
edges with label 1). Conversely, a configuration of (b1)-type with
$f(z)\in2\mathbb Z$ gives
  $$
l=f(z),\quad k=\frac{f(y)-f(z)}2,
  $$
and $f(w)+f(x)-f(y)$ determines the position of $v$ on the path
$P_1(v)$.
\medskip

Now let $l=2t+1$. Then, like the corresponding situation  in Case
1, we first traverse one red edge with label 1, then one green
edge with label $\frac 12$, and then follow the route which is
`parallel' to the route in the case with $l'=l-1$. This gives us
only one change compared with the configuration in the above case,
namely: $f(z):=2t+1=f(w)$ (= twice the number of green edges with
label 1 plus the number of green edges with label $\frac 12$).
Conversely, a configuration of (b1)-type with $f(z)\in2\mathbb
Z+1$ gives
  $$
l=f(z),\quad k=\frac{f(y)-f(z)-1}2,
  $$
and $f(w)+f(x)-f(y)$ determines the position of $v$ on the path
$P_1(v)$.\medskip

Now we have to check that the action of $F_1$ and $F_2$ on the
configurations in $B(H,0)$ agrees with the edges in $K(H,0)$. Let
a red edge emanate from $v$ and end in $v'$. If the label on this
edge is $0$, then we are in Case 1, and hence $f'(w)=f(w)+1$, that
is the case. And if the label is $1$, then we are in Case 2, and
hence $f'(x)=f(x)+1$, that is the case.

Let $v$ and $v'$ be connected by a green edge. Then three cases
are possible.

\medskip
(i) The label of $(v,v')$ is 1. Then the coordinates of $v$ and $v'$ in $K(H,0)$ are $(k+1,l-2,m)$ and
$(k,l,m)$, respectively, and we get $f'(z)=f(z)+2$, and there are no other changes, that is the case.

\medskip
(ii) The label of $(v,v')$ is $\frac 12$. Then the coordinates of $v$ and $v'$ in $K(H,0)$ are $(k,l-1,m)$ and
$(k,l,m)$, but $v$ concerns Case 2 and $v'$ concerns Case 1 (see Picture 5), and we get $f(x)=f(y)+1$,
$f(z)=f(w)$, $f'(y)=f(y)+1=f(x)$ and $f'(z)=f(z)+1$, which is the case. \medskip

(iii) The label of $(v,v')$ is $0$. Then the coordinates of $v$ and $v'$ are $(k-1,l,m+2)$ and $(k,l,m)$, that
yield $f'(y)=f(y)+2$, which is the case.

\medskip
Thus, we established the isomorphism between $K(H,0)$ and
$B(H,0)$. QED
\medskip

{\bf Claim}. $K(0,A)=B(0,\frac A2)$.
\medskip

The configurations corresponding to such intervals are of type (a).

Note that each red path at $k$-th level below
the distinguished diagonal-set (which is located at the $0$-th
level) is constituted from $k$ edges with label 0 and $k$ edges
with label 1.

\medskip
\noindent{\em Case 1}: a vertex $v$ is located above the critical point in the
red path $P_1(v)$ and the ``view from the sky'' of $v$ is located
at the level $k$ and at the position $l$ to the right from the
left "border", $l\le A+1-k$. Thus, each vertex of $K(0,A)$ is
characterized by a triple $(k,l,k'+k)$, $k'\le k$, $0\le l\le A-k$
(if $l=0$, then the vertex is on the left-hand side  border).

On the next picture we illustrate the ``view from the sky'' of the
point with the coordinates $(1,2,1)$ in $K(0,4)$\bigskip

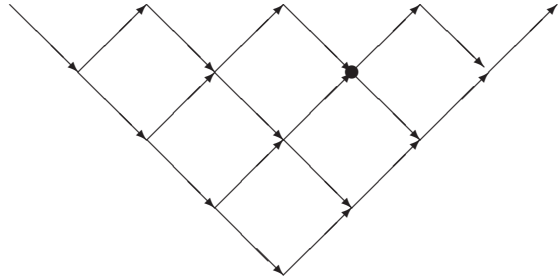
\unitlength=0.6mm \special{em:linewidth 0.4pt}
\linethickness{0.4pt}
\begin{picture}(140.00,65.00)
\put(20.00,60.00){\vector(1,-1){15.00}}
\put(35.00,45.00){\vector(1,-1){15.00}}
\put(50.00,30.00){\vector(1,-1){15.00}}
\put(65.00,15.00){\vector(1,-1){15.00}}
\put(80.00,0.00){\vector(1,1){15.00}}
\put(95.00,15.00){\vector(1,1){15.00}}
\put(110.00,30.00){\vector(1,1){15.00}}
\put(125.00,45.00){\vector(1,1){15.00}}
\put(35.00,45.00){\vector(1,1){15.00}}
\put(50.00,60.00){\vector(1,-1){15.00}}
\put(65.00,45.00){\vector(1,-1){15.00}}
\put(80.00,30.00){\vector(1,-1){15.00}}
\put(50.00,30.00){\vector(1,1){15.00}}
\put(65.00,45.00){\vector(1,1){15.00}}
\put(80.00,60.00){\vector(1,-1){15.00}}
\put(95.00,45.00){\vector(1,-1){15.00}}
\put(65.00,15.00){\vector(1,1){15.00}}
\put(80.00,30.00){\vector(1,1){15.00}}
\put(95.00,45.00){\vector(1,1){15.00}}
\put(110.00,60.00){\vector(1,-1){14.00}}
\put(95.00,45.00){\circle*{3.00}}
\end{picture}

The route from the source to $v$ is of the following form: first
it traverses $k+l$ green edges with label 1, then $k+2l$ red edges,
then $l$ green edges, and finally, $k'$ red edges to enter $v$.

The corresponding admissible configuration is:
$f(x)=f(y)=f(z)=f(w)=0$, $f(a)=2(k+l)$ (= twice the number of
green (or green-blue) edges with label 1); $f(c)=2l$ (= twice the
number of green edges with label 0); $f(b)=k+2l+k'$ (= the number
of red edges on the route (obviously, one holds $f(a)\ge f(b)\ge
f(c)$)). Conversely, we have
  $$
k=\frac{f(a)-f(c)}2,\quad l=\frac{f(c)}2,\quad
k'=f(b)-\frac{f(a)+f(c)}2,
  $$
and $f(b)-f(c)$ is equal to the position of $v$ in $P_1(v)$.
\medskip

\noindent
{\em Case 2}: a vertex $v$ is located below the critical point in the
red path $P_1(v)$. Let $(k,l,k')$ be the coordinates of
$v$. Then our route is as follows: it starts at the source and first
goes to the vertex with the coordinates
$v'=(0,l-k,0)$ if $l\ge k$, and to the vertex $(0,l,0)$ if $l\le
k$, then it goes from $v'$ using green edges with label 1
to the beginning of the string $P_1(v)$, and then makes $k'$ steps
along $P_1(v)$.

As a result of this route, we will traverse $(l-k)$ green edges
with label 1, then $2(l-k)$ red edges, then $(l-k)$ green edges
with label 0 (coming in $v'$), then $k$ green edges with label 1,
and then $k'$ steps to $v$.

So we obtain: $f(a)=2(l-k)+k$ (= twice the number of green
(green-blue) edges with label 1); $f(c)=2(l-k)$ (= the number of
green (green-blue) edges with label 0); $f(b)=k'+2(l-k)$ (= the
number of red edges on the route). Like Case 1, the converse is
true too.
\medskip

A simple verification shows that the action of $F_1$ and $F_2$ in
$B(0,A)$ agrees with the edges in $K(0,\frac A2)$. The Claim is
proven.
\medskip

Finally, the graph $K(H,A)$ contains translated sails of the
types $K(0,A)$ and $K(H,0)$.

The sail $K(H,0)$ attached to the $l$-th distinguished point of
$K(0,A)$, $l\le A$, is described by a copy of the configurations
$K(H,0)$ with the following unimportant modification:
$f(a)=f(b)=f(c)=l$.

As to the sail $K(0,A)$, the situation is a bit tricky.
Specifically, such a sail, attached to the distinguished point of
$K(H,0)$ through the red path $(l,H-l)$, gets the following
modification: modulo the coordinates in the red vertices the
``view from the sky'' of $K_{(l,H-l)}(0,A)$ is, in fact, the same
as of $K(0,A)$, and concerning to the coordinates, all coordinates
of $K_{(l,H-l)}(0,A)$ have to be  change by adding the vector
$(l,H-l)$, that is, by adding $l$ edges from below to the
`original edges' in $K(0,A)$ (with label 0) and $H-l$ edges from
above (with label 1).

The additional vertices in such a sail are managed by the
admissible configurations of type (c1) or (c2), and the
`original' ones are managed by configuration of type (a) with
the only modification $f(x)=f(y)=f(z)=f(w)=l$.

Specifically, if a vertex of $K_{(l,H-l)}(0,A)$ is located below
the 'original vertices' in $K(0,A)$, then $f(x)=f(y)=f(z)=l$, and
$f(w)=l'$, $0\le l'\le l$, and $f(b)=f(c)$  (the (c1)-type
configuration) till we reach original vertices. If a vertex of
$K_{(l,H-l)}(0,A)$ is above the 'original vertices' in $K(0,A)$,
then $f(y)=f(z)=f(w)=l$, and $f(x)=l+l'$ for $0\le l'\le H-l$, and
$f(a)=f(b)$ (the (c2)-type configuration) after the moment we
leave original vertices. A verification of the coincidence of the
crystal operations in $B(H,A/2)$ and $K(H,A)$ in  these additional
cases is by direct calculations.\medskip

Thus, the isomorphism $B(H,A)=K(H,A/2)$ is shown, and this
completes the proof of Theorem 1.

\end{document}